\documentclass[12pt,a4paper]{article}
\usepackage[english]{babel}
\usepackage[utf8]{inputenc}
\usepackage{amsmath}
\usepackage{amsfonts}
\usepackage{amssymb}
\usepackage{graphicx}
\usepackage{listings}
\usepackage{listingsutf8}
\usepackage{tikz}
\usepackage{cite}
\usepackage{algorithm2e}
\usepackage{algpseudocode}
\usepackage{array}
\usepackage{float}
\usepackage{geometry}
\DeclareMathOperator{\trace}{Tr}

\begin{document}

\setcounter{tocdepth}{2}

\nocite{*}
\newcounter{currentcounter}

\newcounter{corollarycounter}

\newcounter{example}
\newcommand\exc{\theexample}

\newenvironment{general}[1]{\refstepcounter{currentcounter} \bigskip \noindent \textbf{#1 \thecurrentcounter :}\itshape}{\normalfont \medskip}

\newenvironment{generalbis}[2]{\refstepcounter{currentcounter} \bigskip \noindent \textbf{#1 \thesubsection.\thecurrentcounter ~#2:}\itshape}{\normalfont \medskip}

\newenvironment{proof}{\noindent \textbf{Proof:} \newline \noindent \hspace*{0.2cm}}{\hspace*{\fill}$\square$ \bigskip}

\newenvironment{proofbis}[1]{\noindent \textbf{#1:} \newline \noindent \hspace*{0.2cm}}{\hspace*{\fill}$\square$ \bigskip}

\newenvironment{famous}[1]{\medskip \noindent \textbf{#1:} \newline \noindent  \itshape}{\normalfont \medskip}

\newenvironment{mycor}{\refstepcounter{corollarycounter} \medskip \noindent \textbf{Corollary \thecorollarycounter :}  \itshape}{\normalfont \medskip}

\newenvironment{court}[1]{\refstepcounter{currentcounter} \smallskip \noindent \textbf{#1 \thesubsection.\thecurrentcounter :} \newline \noindent \itshape}{\normalfont \smallskip}

\newcolumntype{C}{>{$}c<{$}}
\newcolumntype{L}{>{$}l<{$}}
\newcommand\toline{\smallskip \newline}
\newcommand\refp[1]{(\ref{#1})}
\renewcommand\mod{\mathrm{~mod~}}
\newcommand\pregcd{\mathrm{gcd}}
\renewcommand\gcd[2]{\pregcd(#1,#2)}
\newcommand\prelcm{\mathrm{lcm}}
\newcommand\lcm[2]{\prelcm(#1,#2)}
\newcommand\naturaliso{\cong}
\newcommand\iso{\simeq}
\newcommand\cross[1]{#1^{\times}}
\newcommand\crosslong[1]{\cross{(#1)}}
\newcommand\poldegree[1]{\mathrm{deg}(#1)}
\newcommand\cardinal[1]{\# \left(#1\right)}
\newcommand\cardinalshort[1]{\##1}
\renewcommand\det{\mathrm{det}}
\newcommand\indicator[1]{[#1]}

\newcommand\sign{\mathrm{sign}}
\newcommand\signature{\mathrm{sgn}}

\newcommand\bb[1]{\mathbb{#1}}
\newcommand\kk{\bb{K}}
\newcommand\ok{\mathcal{O}_{\kk}}
\newcommand\zz{\bb{Z}}
\newcommand\qq{\bb{Q}}
\newcommand\rr{\bb{R}}
\newcommand\cc{\bb{C}}
\newcommand\ff{\bb{F}}
\newcommand\of{\mathcal{O}_{\ff}}
\newcommand\goth[1]{\mathfrak{#1}}
\newcommand\zsz[1]{\zz/#1\zz}
\newcommand\zszx[1]{(\zsz{#1})^{\times}}
\newcommand\dualsimple[1]{#1^{\vee}}
\newcommand\dual[1]{(#1)^{\vee}}

\newcommand\dx[1]{\mathrm{d}#1}
\newcommand\ddx[1]{\frac{\mathrm{d}}{\dx{#1}}}
\newcommand\dkdxk[2]{\frac{\mathrm{d}^{#2}}{\dx{#1}^{#2}}}
\newcommand\partialx[1]{\partial #1}
\newcommand\partialdx[1]{\frac{\partial}{\partialx{#1}}}
\newcommand\partialkdxk[2]{\frac{\partial^{#2}}{\partialx{#1}^{#2}}}
\newcommand\integral[3]{\int_{#1}{#2}\dx{#3}}
\newcommand\integralsimple[1]{\integrale{G}{#1}{\lambda}}
\newcommand\sprod[2]{\langle #1, #2 \rangle}

\newcommand\limi[1]{\lim_{#1 \to + \infty}}

\newcommand\mn[2]{\mathrm{M}_{#1}(#2)}
\newcommand\mnz[1]{\mn{#1}{\zz}}
\newcommand\sln[2]{\mathrm{SL}_{#1}(#2)}
\newcommand\slnz[1]{\sln{#1}{\zz}}
\newcommand\gln[2]{\mathrm{GL}_{#1}(#2)}
\newcommand\glnz[1]{\gln{#1}{\zz}}
\newcommand\resc{\mathrm{resc}}
\newcommand\com{\mathrm{com}}

\newcommand\bars[1]{\underline{#1}}
\newcommand\taubar{\bars{\tau}}
\newcommand\qbar{\bars{q}}
\newcommand\sigmabar{\bars{\sigma}}
\newcommand\rbar{\bars{r}}
\newcommand\nbar{\bars{n}}
\newcommand\mbar{\bars{m}}
\newcommand\mubar{\bars{\mu}}
\newcommand\xbar{\bars{x}}
\newcommand\Xbar{\bars{X}}
\newcommand\abar{\bars{a}}
\newcommand\bbar{\bars{b}}
\newcommand\alphabar{\bars{\alpha}}
\newcommand\betabar{\bars{\beta}}
\newcommand\omegabar{\bars{\omega}}
\newcommand\taubarsj[1]{\taubar^{-}(#1)}
\newcommand\onebar{\bars{1}}

\newcommand\hh{\bb{H}}
\newenvironment{psmallmatrix}
  {\left(\begin{smallmatrix}}
  {\end{smallmatrix}\right)}
  
\newcommand\declarefunction[5]{#1 := \begin{cases}
\hfill #2 \hfill & \to ~ #3 \\
\hfill #4 \hfill & \to ~ #5 \\
\end{cases}}

\newcommand\sume{\sideset{}{_e}\sum}

\newcommand\tendstowhen[1]{\xrightarrow[#1]{}}

\newcommand\pending{$\bigskip \newline \blacktriangle \blacktriangle \blacktriangle \blacktriangle \blacktriangle \blacktriangle \blacktriangle \blacktriangle \bigskip \newline \hfill$}

\newcommand\pendingref{\textbf{[?]}}

\newcommand\opc[1]{\mathcal{O}^{+,\times}_{#1}}
\newcommand\opck{\opc{\kk}}
\newcommand\opcf{\opc{\goth{f}}}
\newcommand\oonef{\mathcal{O}^{\times}_{\goth{f}}}
\newcommand\eps{\varepsilon}
\newcommand\units{\ok^{\times}}
\newcommand\unitsf{\of^{\times}}
\newcommand\norm[1]{\mathcal{N}(#1)}

\newcommand\ta{\tilde{a}}
\newcommand\tabar{\tilde{\abar}}
\newcommand\talpha{\tilde{\alpha}}
\newcommand\talphabar{\tilde{\alphabar}}
\newcommand\ts{\tilde{s}}
\newcommand\ttt{\tilde{t}}
\newcommand\tlambda{\tilde{\lambda}}
\newcommand\tc{\tilde{c}}
\newcommand\tcone{\tc_1}
\newcommand\tctwo{\tc_2}
\newcommand\te{\tilde{e}}
\newcommand\teone{\te_1}
\newcommand\tetwo{\te_2}
\newcommand\tmu{\tilde{\mu}}
\newcommand\tA{\tilde{A}}

\newcommand\spvgamma{\Gamma_{\goth{f}, \goth{b}, \goth{a}}(\eps; h)}
\newcommand\spvgammasign{\Gamma_{\goth{f}, \goth{b}, \goth{a}}^{\pm}(\eps; h)}
\newcommand\spvgr{G_{r, \goth{f}, \goth{b}, \goth{a}}^{\pm}(u_1, \dots, u_r; h)}
\newcommand\spvgrrho{G_{r, \goth{f}, \goth{b}, \goth{a}}^{\pm}([\eps_{\rho(1)}|\dots|\eps_{\rho(r)}]; h_\rho)}
\newcommand\spvgrcomplete{I_{r,\goth{f}, \goth{b}, \goth{a}}(\eps_1, \dots, \eps_r ; \bars{h}, \bars{\mu}, \bars{\nu})}

\newcommand\spvgammac{\Gamma_{\goth{f}, \goth{b}, \goth{a}}(\eps; h, \sigma_{\cc})}
\newcommand\spvgrc{G_{r, \goth{f}, \goth{b}, \goth{a}}^{\pm}([\eps_{\rho(1)}|\dots|\eps_{\rho(r)}]; h_\rho, \sigma_{\cc})}
\newcommand\spvgrrhoc{G_{r, \goth{f}, \goth{b}, \goth{a}}^{\pm}([\eps_{\rho(1)}|\dots|\eps_{\rho(r)}]; h_\rho)}
\newcommand\spvgrcompletec{I_{r,\goth{f}, \goth{b}, \goth{a}}(\eps_1, \dots, \eps_r ; \bars{h}, \bars{\mu}, \bars{\nu}, \sigma_{\cc})}

\newcommand\ula{u_{L, \goth{a}}}
\newcommand\ulah{u_{L, \goth{a}, \bars{h}}}
    
\newcommand\hombase{\mathrm{Hom}}
\newcommand\myhom[3]{\hombase_{#1}(#2, #3)}    
\newcommand\homlong[3]{\hombase_{#1}\left(#2, #3\right)}   
\newcommand\homlz{\myhom{\zz}{L}{\zz}}
\newcommand\homlprimez{\myhom{\zz}{L'}{\zz}}
\newcommand\homlc{\myhom{\zz}{L}{\cc}}
\newcommand\homlambdaz{\myhom{\zz}{\Lambda}{\zz}}
\newcommand\homlambdac{\myhom{\zz}{\Lambda}{\cc}}
\newcommand\homvq{\myhom{\qq}{V}{\qq}}
\newcommand\zexc{z}
  
\newcommand\zfone{\mathcal{Z}_{\goth{f}}^1}
  
\newcommand\fracpart[1]{\left\{#1\right\}}
\newcommand\entirepart[1]{\left\lfloor#1\right\rfloor}  
\newcommand\shortsetseparator[1][]{#1|}
\newcommand\setseparator[1][]{~\shortsetseparator[#1]~}
\newcommand\coeff{\mathrm{coeff}}
\newcommand\normalbone[1]{\left(\left(#1\right)\right)}
\newcommand\tensor{\otimes}
\newcommand{\at}[2][]{#1|_{#2}}
\newcommand\omitvar[1]{\widehat{#1}}
\newcommand\kronecker{\delta}
\newcommand\plgt{x}
\newcommand\cohomd{\partial}
\newcommand\cohomdx{\partial^{\times}}
\newcommand\compset[1]{#1^c}
\newcommand\signdet{\sign\,\det}

\newcommand\copen{c^{\circ}}
\newcommand\cclosed{c}
\newcommand\cdual{c^{\vee}}
\newcommand\copendual{c^{\vee, \circ}}
\newcommand\dirac{\delta}

\newcommand\myspan{\mathrm{Span}}
\newcommand\prespanconvex[1]{\mathcal{C}(#1)}
\newcommand\spanconvex{\prespanconvex{V}}
\newcommand\prespancones[1]{\mathcal{K}(#1)}
\newcommand\spancones{\prespancones{V}}
\newcommand\spanconesrr{\prespancones{V_{\rr}}}
\newcommand\spanconesk[1]{\mathcal{K}^{#1}(V)}
\newcommand\spanqcones{\mathcal{K}_{\qq}(V)}
\newcommand\spanqconesrr{\mathcal{K}_{\qq}(V_{\rr})}
\newcommand\prespanwedges[1]{\mathcal{L}(#1)}
\newcommand\spanwedges{\prespanwedges{V}}
\newcommand\spanqwedges{\mathcal{L}_{\qq}(V)}
\newcommand\spanqwedgesrr{\mathcal{L}_{\qq}(V_{\rr})}

\newcommand\coefficient[3]{\mathrm{coeff}[#2^{#3}]\left(#1\right)}
\newcommand\badposition{\mathrm{(BP)}}
\newcommand\symbolparagraph{\S}

\newcommand\alphajk[2]{\alpha^{(#1)}_{#2}}
\newcommand\djk[2]{d^{(#1)}_{#2}}
\newcommand\vjkl{v^{(j)}_{k,j'}}
\newcommand\yj[1]{y^{(#1)}}
\newcommand\gammaj[1]{\gamma^{(#1)}}
\newcommand\yjk[2]{y^{(#1)}_{#2}}
\newcommand\ujk[2]{u^{(#1)}_{#2}}

\newcommand\disjointunion{\sqcup}
\newcommand\dkxs{d_k(x, s)}

\newcommand\hexp{h_e}

\newcommand\smoothedgr{G_{n-2, a_1, \dots, a_{n-1}}(v)(w, x, L, L')}
\newcommand\smoothedgromitj{G_{n-2, a_1, \dots, \omitvar{a_j}, \dots, a_{n}}(v)(w, x, L, L')}
\newcommand\smoothedbn{B_{n, a_1, \dots, a_n}(v)(w,x, L, L')}
\newcommand\congruencegrouphnn{H_n(N)}
\newcommand\denomclasscohom{\mathcal{D}(n,N)}

\newcommand\cun[1]{c^1_{#1}}
\newcommand\cdeux[1]{c^2_{#1}}
\newcommand\Cun[1]{C^{1}_{#1}}
\newcommand\Cdeux[1]{C^{2}_{#1}}
\newcommand\fun{f^1}
\newcommand\fdeux{f^2}
\newcommand\mathcalCpm[2]{\mathcal{C}^{#1}_{#2}}
\newcommand\mathcalCplus[1]{\mathcal{C}^{+}_{#1}}
\newcommand\mathcalCmoins[1]{\mathcal{C}^{-}_{#1}}
\newcommand\good{\textit{good }}

\newcommand\preclassgroup{\mathrm{Cl}}
\newcommand\preclassgroupplus{\preclassgroup^{+}}
\newcommand\classgroup[1]{\preclassgroup(#1)}
\newcommand\classgroupplus[1]{\preclassgroupplus(#1)}
\newcommand\classgroupk{\classgroup{\kk}}
\newcommand\classgroupplusk{\classgroupplus{\kk}}
\newcommand\classgroupf{\classgroupplus{\goth{f}}}
\newcommand\wideclassgroupf{\classgroup{\goth{f}}}
\newcommand\classfieldf{\kk^{+}(\goth{f})}

\newcommand\tildeD{\widetilde{\goth{D}}}
\newcommand\zetaf{\zeta_{\goth{f}}}

\begin{center}
\end{center}
\begin{center}
\large \noindent Computations of higher elliptic units in optimal settings
\end{center}
\medskip
\begin{center}
Pierre L. L. Morain\footnotemark[1]\footnotetext[1]{Sorbonne Université and Université Paris Cité, CNRS, INRIA, IMJ-PRG, F-75005 Paris, France. This work is funded by the École polytechnique, Palaiseau, France.}
\end{center}
\medskip
\begin{center}
\noindent \textbf{Abstract:}
\end{center}
In this paper we present a simplified form of a conjecture on the construction of generalised elliptic units above number fields with exactly one complex place. They are conjectural algebraic numbers which are obtained as special values of higher elliptic Gamma functions. These functions form a collection of multivariate meromorphic functions which were studied in the late 1990s and early 2000s in mathematical physics. Our construction extends the scheme of a recent article by Bergeron, Charollois and Garc\'ia where they constructed conjectural elliptic units above complex cubic fields using the elliptic Gamma function. The higher elliptic units we construct are expected to generate specific abelian extensions of the base field where they are evaluated, thus giving a conjectural solution to Hilbert's 12th problem for the number fields with exactly one complex place. We provide several examples to support our conjecture in optimal settings for number fields of degree 3, 4, 5 and 6.

\normalsize

\tableofcontents

\section{Introduction}\label{introduction}

In this paper we are interested in the construction of generalised elliptic units above number fields with exactly one complex place. For imaginary quadratic fields, elliptic units are built in the context of \textit{Complex Multiplication} using a $\theta$ function defined on $\cc \times \hh$ by:
\begin{equation}\label{deftheta}
\theta(z, \tau) = \prod_{n \geq 0} \left(1- e^{-2i\pi z} e^{2i \pi (n+1) \tau}\right) \left(1-e^{2i\pi z}e^{2i\pi n \tau}\right)
\end{equation}
where $\hh$ denotes the upper half-plane. For an imaginary quadratic field $\kk$, evaluations of the form
\begin{equation}\label{ellipticunitsform} 
\frac{\theta\left(\frac{k}{q}, \tau\right)^N}{\theta\left(\frac{N\cdot k}{q}, N\cdot \tau\right)}
\end{equation}
for well-chosen $q, N \in \zz_{\geq 2}$, $k \in \zszx{q}$ and $\tau \in \kk - \qq$ yield $q$-units inside ray class fields of $\kk$ which are called elliptic units (see \cite{Ramachandra}, \cite{Robert}). As a consequence of this construction, Hilbert's 12th problem is solved for imaginary quadratic fields as the $\theta$ function describes the abelian extensions of imaginary quadratic fields. As an example, one may compute the following elliptic units:
\begin{align*}
u_{1} & = \theta\left(\frac{1}{13}, \frac{19+\sqrt{-3}}{182}\right)^7\theta\left(\frac{7}{13}, 7\cdot\frac{19+\sqrt{-3}}{182}\right)^{-1}  = \frac{-3-7\sqrt{-3}}{4} + \sqrt{\frac{-41+17\sqrt{-3}}{8}} \\
u_2 & = \theta\left(\frac{2}{13}, \frac{19+\sqrt{-3}}{182}\right)^7\theta\left(\frac{14}{13}, 7\cdot\frac{19+\sqrt{-3}}{182}\right)^{-1}  = \frac{-3-7\sqrt{-3}}{4} - \sqrt{\frac{-41+17\sqrt{-3}}{8}}
\end{align*}
which are two roots of the polynomial $x^4+3x^3+32x^2+13$, the other two roots being their complex conjugates. Thus, $u_1$ and $u_2$ are $13$-units in a quadratic extension of $\qq(\sqrt{-3})$ which is ramified only above a prime ideal of norm $13$ and they are Galois conjugates over $\qq(\sqrt{-3})$. 

In \cite{thesis}, we generalised this construction of conjectural elliptic units to higher degree number fields with exactly one complex place, extending the construction of conjectural elliptic units above complex cubic fields presented in a recent article by Bergeron, Charollois and Garc\'ia \cite{BCG}. Our construction uses multivariate meromorphic functions enjoying modular transformation properties for higher special linear groups $\slnz{n}$ similar to those of the $\theta$ function for $\slnz{2}$. In the case $n=3$, the function used by Bergeron, Charollois and Garc\'ia is Ruijsenaars' elliptic Gamma function \cite{Ruijsenaars} defined on $\cc\times\hh^2$ by:
\begin{equation}\label{defgamma}\Gamma(z, \tau, \sigma) = \prod_{m,n \,\geq 0}\left(\frac{1-e^{2i\pi((m+1)\tau +(n+1)\sigma-z)}}{1-e^{2i\pi(m\tau+n\sigma+z)}} \right).
\end{equation}
This function was studied at length by Felder and Varchenko \cite{FV} and was shown to belong to a natural hierarchy of multiple elliptic Gamma functions by Nishizawa \cite{Nishizawa}. These are meromorphic functions on $\cc \times \hh^{r+1}$ defined by an infinite product:
\begin{multline}\label{defgr}
G_r(z, \tau_0, \dots, \tau_r) = \prod_{m_0, \dots, m_r \,\geq 0} \left(1 - e^{2i\pi(-z + \sum_{j = 0}^r (m_j+1)\tau_j)}\right) \\ \times \left(1- e^{2i\pi (z + \sum_{j = 0}^r m_j\tau_j)}\right)^{(-1)^r}
\end{multline}
for all $r \in \zz_{\geq 0}$ and we may identify $\theta = G_0$ and $\Gamma = G_1$. These functions are defined by their Fourier product so they are $1$-periodic with respect to each of their arguments, and they enjoy pseudo-periodicity relations in their first argument involving lower degree functions in the hierarchy as:
$$G_r(z + \tau_j, \tau_0, \dots, \tau_r) = G_r(z, \tau_0, \dots, \tau_r) G_{r-1}(z, \tau_0, \dots, \omitvar{\tau_j}, \dots, \tau_r) $$
where $\omitvar{\tau_j}$ means as usual that $\tau_j$ should be omitted. These functions have connections with $q$-polylogs as well as $q$-Pochhammer symbols (see for instance \cite{Nishizawa}) and they may be extended to $\cc \times (\cc - \rr)^{r+1}$ by setting for any $0 \leq j \leq r$:
$$G_r(z, \tau_0, \dots, -\tau_j, \dots, \tau_r) = G_r(z + \tau_j, \tau_0, \dots, \tau_j, \dots, \tau_r)^{-1}.$$
The most compelling property they share is the so-called modular property involving some generalised Bernoulli rational functions in $\qq[z](\omega_1,\dots, \omega_n)$ defined by:
\begin{equation}\label{defbn}
\frac{1}{n!}B_{n,n}(z, \omega_1, \dots, \omega_n) = \coefficient{\frac{e^{zt}}{\prod_{j = 1}^n \left(e^{\omega_jt}-1\right)}}{t}{0}.
\end{equation}
The modular property was proven by Felder and Varchenko \cite{FV} for $n = 3$ and by Narukawa \cite{Narukawa} for general $n$ and may be formulated as:
\begin{equation}\label{modularproperty}
\prod_{j = 1}^n G_{n-2}\left(\frac{z}{\omega_j}, \left(\frac{\omega_k}{\omega_j}\right)_{j \neq k}\right) = \exp\left(-\frac{2i\pi}{n!} B_{n,n}(z, \omega_1, \dots, \omega_n)\right)
\end{equation}
provided that $\omega_k/\omega_j \not\in \rr$ for $1 \leq k \neq j \leq n$. In the light of this relation we built collections of equivariant $(n-1)$-cocycles for specific subgroups of $\slnz{n}$ from the Bernoulli rational functions $B_{n,n}(z, \omega_1, \dots, \omega_n)$ (see \cite{firstpaper}). These $(n-1)$-cocycles are then split by the collection of multiple elliptic Gamma functions and we have shown in \cite{secondpaper} that smoothed versions of these multiple elliptic Gamma functions give rise to equivariant $(n-2)$-cocycles for specific subgroups of $\slnz{n}$. This cohomological interpretation generalises the well-known properties of the $\theta$ function which splits the Dedekind-Rademacher $1$-cocycle. 

Our aim in this paper is to give a simple form of the main conjecture in \cite{thesis} in the simplest cases. We shall explain how to construct, in optimal settings, higher elliptic units above a number field $\kk$ of degree $n$ with exactly one complex place by evaluating
smoothed $G_{n-2}$ functions at points in $\kk$, thus giving an analytic description of the abelian extensions of $\kk$ in the spirit of Hilbert's 12th problem. A more general version of this construction is presented in the author's PhD thesis \cite{thesis} and will be the subject of an upcoming paper \cite{thirdpaper}.

In general, by analogy with formula \refp{ellipticunitsform}, the higher elliptic units above a degree $n$ number field $\kk$ with exactly one complex place should be given by products of the form:
\begin{equation}\label{generalshapeellipticunitsintro}
\prod_{j = 1}^m \frac{G_{n-2}(z_j, \tau_{1, j}, \dots, \tau_{n-1,j})^N}{G_{n-2}(Nz_j, N\tau_{1, j}, \dots, N\tau_{n-1,j})}
\end{equation}
where $N$ is a choice of smoothing index and the evaluation parameters $z_j, \tau_{1, j}, \dots, \tau_{n-1, j}$ should be carefully chosen elements in  $\kk$. Most of the work carried out to formulate a precise conjecture on the algebraic nature of some of these values concerns the choice of these evaluation parameters, as well as the control on the size $m$ of the product. Let us briefly say a few words on these two particular matters. It is not true that all evaluations of the form given in \refp{generalshapeellipticunitsintro} yield algebraic numbers. This contrasts with the situation of imaginary quadratic fields where the evaluation of elliptic units can be carried out for any point $\tau \in \qq(\sqrt{-d})$ by choosing $q \in \zz_{\geq 2}$ correctly. It is therefore important to identify what makes an evaluation of the form \refp{generalshapeellipticunitsintro} ``correct'' and the description of a ``correct'' evaluation procedure has been at the heart of the author's PhD thesis \cite{thesis}. The general construction is quite complicated, but it may be simplified considerably when considering only specific optimal settings for which the product \refp{generalshapeellipticunitsintro} takes the simplest form with the smallest possible size of the product $m$. The evaluation of higher elliptic units in optimal settings (see section \ref{sectionhypotheses}) is the main focus of the present paper. 

 To illustrate our construction we first give a simple example for the complex cubic field $\kk = \qq(e^{2i\pi/3} 10^{1/3})$ where the conjectural elliptic units are given by an evaluation of the elliptic Gamma function. Set $z = e^{2i\pi/3} 10^{1/3}$ and define the two parameters $\tau = -5z^2 - 11z + 5230$ and $\sigma = -2z^2 + z + 2335$. Then the four evaluations:
$$u_k = \frac{\Gamma\left(\frac{k}{5}, \frac{\tau}{1485}, \frac{\sigma}{1485}\right)^{11}}{\Gamma\left(\frac{11k}{5}, \frac{11\tau}{1485}, \frac{11\sigma}{1485}\right)} \approx \begin{cases} -27.5333588... - i\cdot32.7146180...  & \text{ for } k = 1\\ -2.2349933... -i\cdot4.9384566... & \text{ for } k = 2\\-0.0760627...+i\cdot0.1680687... & \text{ for } k = 3 \\-0.0150592...+i\cdot0.0178931... & \text{ for } k = 4\end{cases}$$
may be computed to high precision to be close to the four roots of the palindromic relative polynomial 
$$P = x^4 + (-7z^2 +5z+19)x^3 + (-19z^2 +70z-59)x^2 + (-7z^2 +5z+19)x + 1$$
which defines a $\zsz{4}$ extension $\bb{L}$ of $\kk$ ramified only above the prime ideal of norm $5$ in $\ok$. Alternatively, we may check that they are roots of the following palindromic integral polynomial giving an absolute equation of $\bb{L}$ over $\qq$:
$$P_{abs} = x^{12} + 57x^{11} + 1956x^{10} + 4640x^9 + 35415x^8 - 109818x^7 + 150139x^6-\dots$$
Thus, $u_1, u_2, u_3$ and $u_4$ should be Galois conjugated units in an abelian extension of $\kk$ and they are obtained by evaluating the elliptic Gamma function at specific points in $\kk$. In this example the size of the product in \refp{generalshapeellipticunitsintro} is $m = 1$ which is helpful for computations. 

To illustrate our construction further, we describe a second example for the quartic field $\kk = \qq(z)$ where $z$ is the complex root of the polynomial $x^4 -6x^3-x^2-3x+1$ lying in the upper half-plane. In this setting, the conjectural higher elliptic units are given by a product of $G_2$ functions. Let us consider the six evaluation parameters:
\begin{align*}
&\tau = 5z^3 - 29z^2 - 15z + 87, && \tau' = 2z^3 - 13z^2 + z - 24 = - \rho \\
&\sigma = 6z^3 - 39z^2 + 10z + 47, && \sigma' = 5z^3 - 29z^2 - 15z + 87 = \tau  \\
&\rho = -2z^3 + 13z^2 - z + 24, && \rho' = -2z^3 + 13z^2 + 6z + 143.
\end{align*}
Then the complex number
$$u = \frac{G_2\left(\frac{1}{2}, \frac{\tau}{182}, \frac{\sigma}{182}, \frac{\rho}{182}\right)^{13}}{G_2\left(\frac{13}{2}, \frac{13\tau}{182}, \frac{13\sigma}{182}, \frac{13\rho}{182}\right)} \times \frac{G_2\left(\frac{1}{2}, \frac{\tau'}{182}, \frac{\sigma'}{182}, \frac{\rho'}{182}\right)^{13}}{G_2\left(\frac{13}{2}, \frac{13\tau'}{182}, \frac{13\sigma'}{182}, \frac{13\rho'}{182}\right)}$$
may be computed to high precision and it is found to be close to the root $ \approx 4.1210208... - i\cdot5.0617720...$ of the polynomial 
$$x^8 - 7x^7 + 33x^6 + 49x^5 + 17x^4 + 49x^3 + 33x^2 - 7x + 1$$
which defines an absolute equation of a quadratic extension of $\kk$ over $\qq$. This example is developed in section \ref{sectionsimplestquartic}. In this example the size of the product in \refp{generalshapeellipticunitsintro} is $m = 2$. 

We now give an outline of this paper. In section \ref{sectionarithmetic} we describe optimal settings for the computation of higher elliptic units and we give a simplified form of the general conjecture in \cite{thirdpaper} on the algebraicity of precise evaluations of the $G_{n-2}$ functions at points in a degree $n$ number field with exactly one complex place. In section \ref{sectioncomputing} we explain how we perform these computations to obtain numerical evidence supporting the conjecture in optimal settings, and in section \ref{sectionnumerical} we showcase the conjecture for number fields of degree 3, 4, 5 and 6. \smallskip

\noindent \textbf{Acknowledgments:} 
This paper was written during the preparation of the author's PhD and they would like to thank their advisors Pierre Charollois and Antonin Guilloux for their guidance and their helpful comments.

\section{An algebraicity conjecture on special values of $G_r$ functions}\label{sectionarithmetic}

Let $\kk$ be a number field of degree $n \geq 3$ with exactly one complex place. Let $\sigma_{\cc}$ be one of the two complex embeddings of $\kk$. We shall identify in the rest of this paper $\kk$ with its image $\sigma_{\cc}(\kk)$ in $\cc$. Denote by $\ok$ the ring of integers of $\kk$ and fix a proper ideal $\goth{f}$ of $\ok$. Our goal is to describe units in the narrow ray class field $\classfieldf$ associated to the modulus $\goth{f}\infty_1\dots\infty_{n-2}$ where the $\infty_j$'s are the real places of $\kk$. Recall that $\classfieldf$ is a finite abelian extension of $\kk$ whose Galois group is isomorphic to the narrow ray class group $\classgroupf$ at $\goth{f}$ via the Artin map. This group is defined as the quotient $I(\goth{f})/P^+(\goth{f})$ where $I(\goth{f})$ denotes the group of all fractional ideals in $\kk$ which are coprime to $\goth{f}$ and $P^{+}(\goth{f})$ is the subgroup of $I(\goth{f})$ consisting of principal fractional ideals admitting a totally positive generator congruent to $1$ modulo $\goth{f}$. One may define as usual partial zeta functions attached to classes in $\classgroupf$ by:
$$\zetaf(\goth{c}, s) = \sum_{\substack{\goth{b} \, \triangleleft \, \ok \\ \goth{b} \, \in \, \goth{c}}} \norm{\goth{b}}^{-s}$$
where the sum ranges over integral ideals in the class $\goth{c}$. It is a well-known fact (see for instance \cite{Neukirch}) that these partial zeta functions may be extended to meromorphic functions on $\cc -\{1\}$ and that their order of vanishing at $s = 0$ does not depend on the class $\goth{c}$.

Throughout this paper we will consider a class $\goth{c}$ in $\classgroupf$ and represent it by an integral ideal of the form $k\cdot \goth{b}$ where $k$ is a rational integer and $\goth{b}$ is an integral ideal in $\kk$. We shall also need an auxiliary integral ideal $\goth{a}$ of norm $N$ such that $N$ is coprime to both $\goth{f}$ and $\goth{b}$. This ideal $\goth{a}$ will be used to perform a standard smoothing operation inspired by \cite{CassouNogues}, \cite{DasguptaShintani} and \cite{CD}. For the purpose of this paper, we shall assume that $N$ is a prime number. 

The construction of higher elliptic units for $\kk$ is a procedure which to the data of $\goth{f}$, $\goth{c}$ and $\goth{a}$ associates a complex number obtained by evaluating a product of $\goth{a}$-smoothed $G_{n-2}$ functions at points in $\kk$. The general conjecture in \cite{thesis} states that this complex number is a unit inside the class field $\classfieldf$ (which we identify with one of its complex embeddings above $\sigma_{\cc}(\kk)$). Our procedure is inspired by the construction of conjectural elliptic units above complex cubic fields by Bergeron, Charollois and Garc\'ia (see \cite{BCG}). Their construction should in fact produce smoothed Stark units and this is supported by their Kronecker limit formula relating the values of derivatives of partial zeta functions at $s=0$ to values of the elliptic Gamma functions. 

As mentioned in the introduction, higher elliptic units should be given in general by finite products of the form:
\begin{equation}\label{generalshapeellipticunits}
\prod_{j = 1}^m \frac{G_{n-2}(z_j, \tau_{1,j}, \dots, \tau_{n-1,j})^N}{G_{n-2}(Nz_j, N\tau_{1, j}, \dots, N\tau_{n-1,j})}
\end{equation}
where the evaluation parameters $z_j, \tau_{1, j}, \dots, \tau_{n-1, j}$ should be carefully chosen elements in $\kk$. Let us give a bit of context as to why we expect such evaluations to be arithmetically interesting. The idea developed in \cite{thesis} is to evaluate an $(n-2)$-cocycle built from the $G_{n-2}$ functions against an $(n-2)$-cycle associated with a choice of fundamental units for the group $\opcf$ of totally positive units of $\ok$ which are congruent to $1$ modulo $\goth{f}$. This cycle is an analogue for number fields with exactly one complex place of Sczech's $(n-1)$-cycle (see \cite{Sczech}) for totally real number fields. Indeed, $\opcf$ is a multiplicative $\zz$-module of rank $n-2$ by Dirichlet's unit theorem. Furthermore, as the only roots of unity belonging to $\kk$ are $-1$ and $+1$ and because $-1 \not\in \opcf$, $\opcf$ is a free $\zz$-module. Thus, if we fix a set $\eps_1, \dots, \eps_{n-2}$ of fundamental units for $\opcf$, that is, $\opcf = \eps_1^{\zz} \dots \eps_{n-2}^{\zz}$, then we may consider the corresponding $(n-2)$-cycle $\Upsilon \in H_{n-2}(\opcf, \zz)$ given by:
\begin{equation}\label{defcycle}
\Upsilon = \sum_{\rho \in \goth{S}_{n-2}} \signature(\rho)[\eps_{\rho(1)} | \dots | \eps_{\rho(n-2)}]
\end{equation}
where $\goth{S}_{n-2}$ is the symmetric group on $n-2$ elements and for general units $u_1, \dots, u_{n-2}$:
$$[u_1| \dots | u_{n-2}] = (1, u_1, u_1u_2, \dots, \prod_{j = 1}^k u_j, \dots, \prod_{j = 1}^{n-2}u_j).$$
The cycle $\Upsilon$ is actually independent of the choice of fundamental units for $\opcf$, but the rest of our construction is not. The expected shape for the higher elliptic units is then more specifically given by:
\begin{equation}\label{expectedshapedeux}
 \prod_{\rho \in \goth{S}_{n-2}} \prod_{j = 1}^{t_{\rho}}\frac{G_{n-2}\left(z_{\rho, j}, \tau_{1, \rho}, \dots, \tau_{n-1, \rho}\right)^N}{G_{n-2}\left(Nz_{\rho,j}, N\tau_{1, \rho}, \dots, N\tau_{n-1, \rho}\right)}\end{equation}
where the product contains 
$$m = \sum_{\rho \in \goth{S}_{n-2}} t_{\rho} \geq (n-2)!$$
terms, each corresponding in an explicit way to a summand $[\eps_{\rho(1)}|, \dots, | \eps_{\rho(n-2)}]$ of the cycle $\Upsilon$ (see section \ref{sectionsimpleconj}). The optimal case occurs when $m = (n-2)!$ as this gives us the shortest possible product associated to this cycle. Explicitly:

\begin{center}
\noindent \begin{tabular}{| c | C | C | C | C | C |}
\hline
degree $n$ & 3 & 4 & 5 & 6 & 7+ \\
\hline 
minimal number $m$ of terms & 1 & 2 & 6 & 24 & (n-2)! \geq 120 \\
\hline
\end{tabular} 
\end{center}

In general, the conjecture on the algebraic nature of some of these evaluations is quite complicated to state, so we present a simplified form in this paper (see \cite{thesis} for a general statement) under some simplifying hypothesis which we now explain. 

\subsection{Optimal settings}\label{sectionhypotheses}

To state our simplifying hypotheses we introduce notations from \cite{thesis}. Recall that we have fixed a number field $\kk$ (of degree $n \geq 3$) with exactly one complex place, a proper integral ideal $\goth{f}$ and a set $\eps_1, \dots, \eps_{n-2}$ of fundamental units for the group $\opcf$. For any $\rho \in \goth{S}_{n-2}$ and any $1 \leq j \leq n-2$ let us set $u_{j, \rho} = \prod_{i = 1}^{j} \eps_{\rho(i)}$. We may define a few important objects attached to each of the unit systems $u_{1, \rho}, \dots, u_{n-2, \rho}$ following \cite{thesis}. First, we consider the linear form:
$$\tilde{f}_{\rho} := \begin{cases} \kk &\to \qq \\ x & \mapsto \det(1, u_{1, \rho}, \dots, u_{n-2, \rho}, x)\end{cases} $$
where the determinant is taken relative to a fixed $\zz$-basis of $\ok$. The image of $\ok$ by the linear form $\tilde{f}_{\rho}$ is a set of the form $\tlambda_{\rho} \,\zz$ where $\tlambda_{\rho} \in \zz_{\geq 0}$ is what we call the ``content'' of the unit system $u_{1, \rho}, \dots, u_{n-2, \rho}$. Assume for now that this content is not zero. Then, the linear form $\ta_{\rho} = \tilde{f}_{\rho}/\tlambda_{\rho}$ maps $\ok$ to $\zz$. We may associate to this linear form $\ta_{\rho}$ a generalised different ideal $\tildeD_{\rho}$ defined by:
\begin{equation}\label{defoverflowideal}\tildeD_{\rho}^{-1} = \{x \in \kk \setseparator \forall\, y \in \ok, \,\ta_{\rho}(x \cdot y) \in \zz \}.
\end{equation}
The analogy with the usual different ideal of a number field $\goth{d}$ which satisfies
$$\goth{d}^{-1} =\{x \in \kk \setseparator \forall\, y \in \ok, \,\trace(x \cdot y) \in \zz \}$$
is clear. As a matter of fact, using the isomorphism
$$\begin{cases} \kk &\to \myhom{\qq}{\kk}{\qq} \\ \xi &\mapsto (y \mapsto \trace(\xi \cdot y)) \end{cases} $$
we may define $\xi_{\rho} \in \kk$ to be the unique element satisfying $\ta_{\rho} = (y \mapsto \trace(\xi_{\rho} \cdot y))$. In this case, the integral ideal $\tildeD_{\rho}$ is given by $\tildeD_{\rho} = \xi_{\rho} \cdot \goth{d}$. In particular, the class of the ideal $\tildeD_{\rho}$ in the (wide) class group of $\kk$ is independent of the permutation $\rho$. The integral ideal $\tildeD_{\rho}$ is called the ``overflow ideal'' (of the unit system $u_{1, \rho}, \dots, u_{n-2, \rho}$) and the positive integer $\ttt_{\rho}$ satisfying $\tildeD_{\rho} \cap \zz = \ttt_{\rho} \, \zz$ is called the ``overflow'' (of the unit system). Let us now regroup the hypotheses made in \cite{thesis} as well as simplifying hypotheses we impose in the present paper to achieve a minimal formulation of the conjecture in optimal cases. Let us start with the general hypotheses of \cite{thesis}:
\begin{enumerate}
\item Assume $\goth{f} \neq \ok$. As it is already the case for classical elliptic units above imaginary quadratic fields, the class field modulus $\goth{f} = \ok$ behaves very differently and requires specific adaptations which we won't discuss here.
\item Assume that $\ok/\goth{f}$ is a cyclic group (this is [\hspace{1sp}\cite{thesis}, (H5)]). This is a simplifying hypothesis which plays a role in many technical aspects of the construction. 

\item Assume that for any $\rho \in \goth{S}_{n-2}$, the content $\tlambda_{\rho}$ of the unit system $u_{1, \rho}, \dots, u_{n-2, \rho}$ is exactly $\tlambda_{\rho} = 1$ (this is [\hspace{1sp}\cite{thesis}, (H2) and (H3)]). This hypothesis is here to simplify the constructions considerably and allows us to consider higher elliptic units given by a product containing the minimal $(n-2)!$ number of terms. Such an assumption is quite restrictive. 
\item Assume that for any $\rho \in \goth{S}_{n-2}$, the overflow ideal $\tildeD_{\rho}$ is coprime to $\norm{\goth{f}}$ (this is [\hspace{1sp}\cite{thesis}, (H4)]). This also simplifies the situation for the computations, and this hypothesis is satisfied quite often.
\end{enumerate}

The main conjecture in \cite{thesis} is expressed under these simplifying hypotheses which we specifically introduced for that matter. Extensive numerical computations of conjectured forms for the higher elliptic units were performed in order to identify the \textit{contents} and \textit{overflow ideals} as key parameters for our evaluations. To present a simpler version of it we shall make three additional assumptions:

\begin{enumerate}
\setcounter{enumi}{4}
\item Assume that $\opck = \opcf$, that is, all totally positive units in $\ok$ are congruent to $1 \mod \goth{f}$. This is the ``optimal setting'' hypothesis \textit{per se}. This is a very restrictive hypothesis as for a given number field $\kk$, only finitely many such ideals exist (see lemma \ref{lemmafiniteJ}).
\item Assume that $\kk$ does not contain any subfield (equivalently, the Galois group over $\qq$ of the Galois closure of $\kk$ is $\goth{S}_n$). In the case where $\kk$ contains a subfield (which is automatically totally real) there are technical difficulties to overcome in the construction. Since we are interested in presenting the simplest possible examples, we shall not consider this case.
\item Assume that the class of the overflow ideal $\tildeD_{\rho}$ in the narrow Hilbert class group of $\kk$ is independent of the permutation $\rho$. This hypothesis is expected to be true generally and it allows us to use the alternative form of the conjecture in \cite{thesis}.
\end{enumerate}

One of the benefits of the ``optimal setting'' hypotheses is that we have a simplified exact sequence in class field theory to work with. Indeed, consider the well-known exact sequence given by: 
\begin{equation}\label{exactseqclgpplus}
1 \to \opcf \to \opck \to (\ok/\goth{f})^{\times} \to \classgroupf \to \classgroupplusk \to 1
\end{equation}
where $\classgroupplusk$ is the narrow Hilbert class group of $\kk$. Since we have assumed that $\opck = \opcf$ and that $(\ok/\goth{f})^{\times} \simeq \zszx{q}$, the exact sequence \refp{exactseqclgpplus} becomes:
\begin{equation}\label{exactseqsimple}
1 \to \zszx{q} \to \classgroupf \to \classgroupplusk \to 1.
\end{equation}
Note that for a given number field $\kk$, there are only finitely many integral ideals $\goth{f}$ satisfying $\opcf = \opck$ as a consequence of the following lemma:

\begin{general}{Lemma}\label{lemmafiniteJ}
Let $\kk$ be a number field and let $\eps_1, \dots, \eps_m$ be generators of the abelian group $\opck$. Then the ideal 
\begin{equation}\label{defjk}J(\kk) = \sum_{j = 1}^{m} (\eps_j -1) \ok = \sum_{\eps \in \opck} (\eps -1) \ok\end{equation}
is independent of the choice of generators $\eps_1, \dots, \eps_m$. An integral ideal $\goth{f}$ satisfies $\opcf = \opck$ if and only if $\goth{f}$ is a divisor of $J(\kk)$.
\end{general}

\begin{proof}
Let us denote temporarily by $J(\eps_1, \dots, \eps_m)$ the integral ideal $\sum_{j = 1}^{m} (\eps_j -1) \ok $ associated to any finite set of generators $\eps_1, \dots, \eps_m$ for the group $\opck$. Let us first prove that $J(\eps_1, \dots, \eps_m)$ divides the integral ideal $(\eps-1)\ok$ for any $\eps \in \opck$. Indeed, let $\eps \in \opck$ be any totally positive unit. There are integers $k_1, \dots, k_m$ such that $\eps = \prod_{j = 1}^m \eps_j^{k_j}$ and therefore the ideal
$$(\eps -1)\ok = \prod_{j = 1}^m (\eps_j -1)\left(\sum_{l_j = 0}^{k_j-1} \eps_j^{l_j}\right) \ok = J(\eps_1, \dots, \eps_m)\left[\prod_{j = 1}^m\left(\sum_{l_j = 0}^{k_j-1} \eps_j^{l_j}\right) \ok\right]  $$
is divisible by $J(\eps_1, \dots, \eps_m)$. In particular, for any other finite set of generators $(\eps'_i)_{i \in I}$, 
$$J(\eps_1, \dots, \eps_m) \subset \sum_{i \in I} (\eps'_i -1)\ok = J(\eps'_i, i \in I).$$
Since the set of generators $\eps_1, \dots, \eps_m$ was arbitrary, we may repeat this part of the proof with the set $(\eps'_i)_{i \in I}$ and obtain the converse inclusion $J(\eps'_i, i \in I) \subset J(\eps_1, \dots, \eps_m)$. Thus the ideal $J(\kk) = J(\eps_1, \dots, \eps_m)$ is independent from the choice of generators for $\opck$. Since it divides any integral ideal of the form $(\eps -1)\ok$ for $\eps \in \opck$ we may write it informally as 
$$J(\kk) = \sum_{\eps \in \opck} (\eps-1)\ok.$$
Let us now suppose that $\goth{f}$ is an integral ideal satisfying $\opcf = \opck$. Then $\goth{f}$ divides all the integral ideals $(\eps_j -1)\ok$ for $1 \leq j \leq m$ and therefore $\goth{f}$ divides the integral ideal $J(\kk) = \sum_{j = 1}^m (\eps_j-1)\ok$. Conversely, if $\goth{f}$ divides $J(\kk)$ then $\goth{f}$ divides all integral ideals $(\eps -1) \ok$ for $\eps \in \opck$ and all units $\eps \in \opck$ belong to $\opcf$. This completes the proof.
\end{proof}

In our setting, $\opck$ contains no torsion, so we may fix a set $\eps_1, \dots, \eps_{n-2}$ of fundamental units for $\opck$, so that the ideal $J(\kk)$ is given by $J(\kk) = \sum_{j = 1}^{n-2} (\eps_j -1)\ok$. 

\subsection{Conjectural elliptic units in optimal settings}\label{sectionsimpleconj}

To state the simplified form of the main conjecture in \cite{thesis} for optimal settings we now describe the evaluation parameters $z_{\rho, j}, \tau_{l, \rho}$ appearing in the evaluation of the higher elliptic units described in \refp{expectedshapedeux}. Recall that we have fixed the modulus $\goth{f}$ of norm $q > 1$, the class $\goth{c}$ represented by the integral ideal $k \cdot \goth{b}$ with $k \in \zz_{\geq 1}$ and $\goth{b}$ an integral ideal, as well as the smoothing ideal $\goth{a}$ of prime norm $N$. Recall that $\goth{f}$ is coprime to $k \cdot \goth{b}$ and to $N$. Recall that we have also fixed a set of fundamental units for $\opcf = \opck$. To describe precisely the evaluation parameters we need to introduce the notion of ``helper ideals'' and the notion of ``compatible sets of base points'' from \cite{thesis}. 

\begin{general}{Definition}
A helper ideal is any integral ideal which is either $\ok$ or of the form $p\ok \times \goth{p}^{-1}$ where $p$ is a rational prime (coprime to $q \cdot N$) and $\goth{p}$ is a prime ideal of norm $p$.
\end{general}

\begin{general}{Definition}\label{defcompatible}
Let $(h_{\rho})_{\rho}$ be a family of totally positive elements in $\goth{f}\goth{b}^{-1}$ such that for any $\rho \in \goth{S}_{n-2}$, $h_{\rho}/q \equiv 1 \mod \goth{f}\goth{b}^{-1}$. The family $(h_{\rho})_{\rho}$ is said to be a compatible set of base points if for any permutation $\rho \in \goth{S}_{n-2}$:
$$h_{\rho} \ok = m_{\rho}\cdot \frac{qN}{\goth{a}\goth{b}}\cdot \tildeD_{\rho}\cdot \goth{H} $$
where $m_{\rho}$ is a positive integer coprime to $q$, the integral ideal $\tildeD_{\rho}$ is the overflow ideal defined in \refp{defoverflowideal} and $\goth{H}$ is a common helper ideal independent from the permutation $\rho$.
\end{general}

These two notions were introduced for the purpose of the construction of higher elliptic units and they play a crucial role in the explicit definition of the evaluation parameters. Let us now prove that such families exist under the hypotheses described in section \ref{sectionhypotheses}.

\begin{general}{Proposition}
Under the hypotheses of section \ref{sectionhypotheses}, there is a compatible set of base points $(h_{\rho})_{\rho}$ for the data $\goth{f}, \goth{b}, \goth{a}, \eps_1, \dots, \eps_{n-2}$.
\end{general}

\begin{proof}
First, we note that any class in the narrow Hilbert class group $\classgroupplusk$ of $\kk$ contains infinitely many helper ideals. Indeed, by Cebotarev's density theorem, each class contains infinitely many prime ideals of prime norm, and, up to removing a finite amount of them, each class contains infinitely many prime ideals of prime norm coprime to $q \cdot N$. Consider $\goth{C}$ a class in the narrow ray class group of $\kk$. For any prime ideal $\goth{p}$ of prime norm $p$ in the class $\goth{C}^{-1}$ such that $(p, q\cdot N) = 1$, the ideal $p\ok \times \goth{p}^{-1}$ is a helper ideal in $\goth{C}$. Hence the class $\goth{C}$ contains infinitely many helper ideals.

Now, by assumption, the class of the fractional ideal $\frac{N}{\goth{a}\goth{b}} \cdot\tildeD_{\rho}$ in $\classgroupplusk$ is independent of $\rho$. Thus there is a helper ideal $\goth{H}$ such that for each $\rho$ the fractional ideal $\frac{N}{\goth{a}\goth{b}} \cdot \tildeD_{\rho} \cdot \goth{H}$ belongs to the trivial class of $\classgroupplusk$. Using the exact sequence \refp{exactseqsimple} we get that the class of the fractional ideal $\frac{N}{\goth{a}\goth{b}} \cdot \tildeD_{\rho} \cdot \goth{H}$ in $\classgroupf$ is the class of some ideal $k_{\rho} \ok$ for some positive integer $k_{\rho}$ coprime to $q$. Let us fix for each $\rho$ a positive integer $m_{\rho}$ such that $m_{\rho} \cdot k_{\rho} \equiv 1 \mod q$. Then each fractional ideal $m_{\rho}\cdot \frac{N}{\goth{a}\goth{b}} \cdot\tildeD_{\rho} \cdot\goth{H} $ belongs to the trivial class of the narrow ray class group at $\goth{f}$ and is therefore generated by some totally positive element $g_{\rho}$ satisfying $g_{\rho} \equiv 1 \mod \goth{f}$. Setting $h_{\rho} = q\cdot g_{\rho}$ completes the proof.
\end{proof}

From now on, we consider a fixed compatible set of base points $(h_{\rho})_{\rho}$. Attached to this data we may consider the linear form
$$f_{\rho} := \begin{cases} \kk &\to \qq \\ x & \mapsto \det(h_{\rho}, u_{1, \rho}h_{\rho}, \dots, u_{n-2, \rho}h_{\rho}, x)\end{cases} $$
as well as the attached rational polyhedral cone
$$C_{\rho} := \{ y \in \kk \setseparator f_{\rho}(y) \geq 0, f_{\rho}(u_{1, \rho}\cdot y) \geq 0, \dots, f_{\rho}(u_{n-2, \rho}\cdot y) \geq 0\}.$$
The intersection of this cone with $\goth{f}(k \cdot \goth{b})^{-1}$ is given by:
$$C_{\rho} \cap \goth{f}(k \cdot \goth{b})^{-1} = \frac{1}{k} \left(\zz \frac{h_{\rho}}{m_{\rho}} + \zz_{\geq 0} \alpha_{1, \rho} + \dots +\zz_{\geq 0} \alpha_{n-1, \rho} \right) $$
for some elements $\alpha_{1, \rho}, \dots, \alpha_{n-1, \rho} \in \goth{f}(k \cdot \goth{b})^{-1}$ (see \cite{thesis}). In particular, the cone $C_{\rho} \cap \goth{f}(k \cdot \goth{b})^{-1}$ is unimodular and the integer $t_{\rho}$ appearing in \refp{expectedshapedeux} is equal to $1$ under our working hypotheses. The conjectural elliptic units are given by a formula of the shape
\begin{equation}\label{expectedshapetrois}
 \prod_{\rho \in \goth{S}_{n-2}}\frac{G_{n-2}\left(z_{\rho}, \tau_{1, \rho}, \dots, \tau_{n-1, \rho}\right)^N}{G_{n-2}\left(Nz_{\rho}, N\tau_{1, \rho}, \dots, N\tau_{n-1, \rho}\right)}\end{equation}
where the evaluation parameters are explicitly given by:
\begin{equation}\label{defzparam} 
z_{\rho} = \epsilon_{\rho} \frac{k\cdot m_{\rho}}{q} \end{equation}
\begin{equation}\label{deftauparam}
\tau_{j, \rho} = \epsilon'_{\rho} \frac{m_{\rho} \cdot \alpha_{j, \rho}}{h_{\rho}} \end{equation}
for some signs $\epsilon_{\rho}, \epsilon'_{\rho} \in \{\pm 1\}$. We may now formulate the simple form of the conjecture in optimal settings.

\begin{generalbis}{Conjecture}{[Simplified form of [\hspace{1sp}\cite{thesis}, Conjecture III.38]]}\label{conjecture}
\newline\noindent Let $\kk$ be a number field of degree $n \geq 3$ with exactly one complex place. Fix an integral ideal $\goth{f} \neq (1)$ and set $q = \norm{\goth{f}}$. Assume that $(\kk, \goth{f})$ is an optimal setting. Fix a class $\goth{c}$ in $\classgroupf$ and a representation of this class by an integral ideal of the form $k \cdot \goth{b}$ for some positive integer $k$ and some integral ideal $\goth{b}$. Fix a smoothing ideal $\goth{a}$ of prime norm $N > n+1$ such that $N$ is coprime to both $\goth{f}$ and $\goth{b}$. Fix a set $\eps_1, \dots, \eps_{n-2}$ of fundamental units for the group $\opcf = \opck$. Assume that the hypotheses of section \ref{sectionhypotheses} are satisfied. Fix a compatible set of base points $(h_{\rho})_{\rho}$. Define the evaluation parameters $z_{\rho}, \tau_{1, \rho}, \dots, \tau_{n-1, \rho}$ as in \refp{defzparam} and \refp{deftauparam} for some orientation signs $\epsilon_{\rho}, \epsilon'_{\rho} \in \{\pm 1\}$. Define the complex number 
\begin{equation}\label{defukb}
 u_{\goth{f}, k, \goth{b}, \goth{a}}(\eps_1, \dots, \eps_{n-2}, (h_{\rho})_{\rho}, (\epsilon_{\rho}, \epsilon'_{\rho})_{\rho}) =  \prod_{\rho \in \goth{S}_{n-2}}\frac{G_{n-2}\left(z_{\rho}, \tau_{1, \rho}, \dots, \tau_{n-1, \rho}\right)^N}{G_{n-2}\left(Nz_{\rho}, N\tau_{1, \rho}, \dots, N\tau_{n-1, \rho}\right)}.
\end{equation}
Then there are orientation signs $\epsilon_{\rho}, \epsilon'_{\rho}$ such that the complex number 
$$u_{k, \goth{b}} = u_{\goth{f}, k, \goth{b}, \goth{a}}(\eps_1, \dots, \eps_{n-2}, (h_{\rho})_{\rho}, (\epsilon_{\rho}, \epsilon'_{\rho})_{\rho})$$
 is an algebraic unit inside $\classfieldf$ which is independent of the choice of fundamental units and compatible set of base points. In addition:
\begin{itemize}
\item The family $(u_{k, \goth{b}})_{k \in \zszx{q}}$ obtained when varying $k$ consists of units in $\classfieldf$ which are Galois conjugates over the narrow Hilbert class field of $\kk$.
\item The unit $u_{\goth{c}} = u_{k, \goth{b}}$ is independent of the choice of representation of $\goth{c}$ by the integral ideal $k \cdot \goth{b}$.
\item The family $(u_{\goth{c}})_{\goth{c} \in \classgroupf}$ obtained when varying the class $\goth{c}$ over $\classgroupf$ consists of Galois conjugated units over $\kk$. The explicit reciprocity law is given by $\sigma_{\goth{c}}(u_{\goth{c}'}) = u_{\goth{c}\goth{c}'}$ where $\goth{c} \mapsto \sigma_{\goth{c}}$ is the Artin map.
\item The unit $u_{\goth{c}}$ satisfies a Kronecker limit formula of the form:
\begin{equation}\label{klfconjecture}
N\zeta'_{\goth{f}}(\goth{c}, 0) - \zeta'_{\goth{f}}(\goth{a}\goth{c}, 0) = \log|u_{\goth{c}}|^2.
\end{equation}
\end{itemize}
\end{generalbis}

\clearpage

\textbf{Remarks on the conjecture:}
\begin{enumerate}
\item This conjecture contains two main statements: an algebraicity property for $u_{k, \goth{b}}$ with applications to Hilbert's 12th problem and a so-called Kronecker limit formula by analogy with the famous second limit formula due to Kronecker. In the case of an imaginary quadratic field, the elliptic units given in \refp{ellipticunitsform} are known to be $q$-units in $\classfieldf$ (see \cite{Ramachandra}, \cite{Robert}) and for a particular choice of $\tau \in \qq(\sqrt{-d})$ they satisfy Kronecker's second limit formula. In the case of complex cubic fields, the Kronecker limit formula was proven by Bergeron, Charollois and Garc\'ia in \cite{BCG}, but the algebraicity statement (which is [\hspace{1sp}\cite{BCG}, Conjecture]) remains open.
\item The construction of these higher elliptic units is inspired by the rank one abelian Stark conjecture for number fields with exactly one complex place, and the units $u_{\goth{c}}$ may be viewed as conjugates of a smoothed version of the conjectural Stark unit. We refer to \cite{rankoneStark} for a more complete presentation of the rank one abelian Stark conjectures. We also refer to \cite{roblot}, \cite{CharolloisDarmon}, \cite{DasguptaShintani}, for computations of conjectural Stark units in various contexts.
\item When the partial zeta functions all vanish at $s = 0$ with order at least $2$, it is expected that the complex number $u_{\goth{c}}$ is a root of unity. The statement is still non-trivial in this case, but examples where this isn't the case are more interesting as the higher elliptic units may be algebraic numbers for which no other parametrisation is known yet.
\item The general case of the conjecture (see \cite{thesis}) without some of the assumptions on $\goth{f}$ is more complicated to state. Indeed, the precise computation of the $u_{\goth{c}}$'s depends on the properties of the cycle $\Upsilon$ which depends on a choice of fundamental units. A poor choice of fundamental units may lead to very lengthy computations and overall to complications in the choice of the base points and evaluation parameters. 
\item The condition that the smoothing index $N$ is greater than $n+1$ can be removed, in which case the Kronecker limit formula should still hold, but the unit $u_{\goth{c}}$ should lie in a cyclotomic extension of $\classfieldf$. 
\item This conjecture was formulated based upon a great number of numerical examples in degrees $n = 3, 4, 5$ (see section \ref{sectionnumerical}) and a few examples in degree $n = 6$. Despite our best efforts, we know not of any number field of degree $7$ with exactly one complex place where the simplifying hypothese described in section \ref{sectionhypotheses} for all $\rho \in \goth{S}_{5}$ are satisfied (and regardless, of any number field of degree $7$ where the computations to test the conjecture take a reasonable amount of time).  
\item This conjecture is expressed up to the signs $\epsilon_{\rho}, \epsilon'_{\rho}$ which we determine by checking the Kronecker limit formula \refp{klfconjecture}. We expect that these signs may be expressed explicitly in terms of related signs introduced by Espinoza \cite{Espinoza} for their explicit Shintani cones for number fields with exactly one complex place.
\end{enumerate}

In the next section, we will explain briefly how we obtain numerical evidence to support the conjecture and what difficulties we need to overcome to perform the computations. We shall then present examples of computations for number fields of degree 3, 4, 5 and 6.

\section{Computational aspects}\label{sectioncomputing}

\subsection{Computing the $G_r$ functions}\label{sectioncomputinggr}

Let us now explain how we compute the $G_r$ functions, practically speaking, to verify the conjecture on numerical examples. Indeed, the definition of the $G_r$ functions by a multi-index infinite product given by formula \refp{defgr} is not suited for computations. In the case $n = 2$, Jacobi's triple product formula gives a beautiful expression of the $\theta$ function in terms of an infinite sum with converging rate in $q^{n^2/2}$:
\begin{equation}\label{jacobitriple}
\theta(z, \tau) = \frac{q^{1/24}}{\eta(\tau)} \sum_{n \in \zz} x^n (-1)^n q^{n(n-1)/2}
\end{equation}
where $x = \exp(2i\pi z)$, $q = \exp(2i\pi\tau)$ and $\eta(\tau)$ is Dedekind's $\eta$ function. This makes the computation of $\theta$ very fast when $\tau$ is not too close to the real axis, and standard techniques in the study of modular functions allow one to use the modularity of $\theta$ to reduce the computation of $\theta(z, \tau)$ to the computation of $\theta(z', \tau')$ where $\Im(\tau') \geq 1/2$. Unfortunately, there is no clear generalisation of Jacobi's triple product formula for higher degree $G_r$ functions and we need other techniques to compute the $G_r$ functions efficiently. We shall make use of [\hspace{1sp}\cite{Nishizawa}, Proposition 3.6] which we write as:
\begin{equation}\label{sumcomputations}
G_r(z, \tau_0, \dots, \tau_r) = \exp\left(-\sum_{j \geq 1} \frac{1}{j} \frac{q_0^j\dots q_r^jx^{-j} +(-1)^r x^{j}}{\prod_{k = 0}^r (1 - q_k^j)}\right)
\end{equation}
where $x = \exp(2i\pi z)$ and $q_k = \exp(2i\pi \tau_k)$ for $0 \leq k \leq r$.
This formula is only valid for $\tau_0, \dots, \tau_r \in \hh^{r+1}$ and $0 < \Im(z) < \sum_{k = 0}^r \Im(\tau_k)$. We call this domain the center strip. It is remarkable that the complexity in the computation of the right-hand side in \refp{sumcomputations} is linear in $r$, whereas usual computations in algebraic number theory tend to be exponential in the degree $n$ of the number field. To compute the value $G_r(z, \tau_0, \dots, \tau_r)$ in general, we only need to make sure that we can use the properties of the $G_r$ functions to reach this domain. This requires two steps in general: a reorientation step and a translation step. 

The reorientation step consists of reducing the computation of any value $G_r(z, \tau_0, \dots, \tau_r)$ with parameters $\tau_0, \dots, \tau_r \in \cc - \rr$ to the computation of some $G_r(z', \tau_0', \dots, \tau_r')$ with parameters $\tau_0', \dots, \tau_r'$ in the upper half-plane. To achieve this, we use repeatedly the inversion property:
\begin{equation}\label{inversionpropertyalgo}
G_r(z, \tau_0, \dots, \tau_{k-1}, -\tau_k, \tau_{k+1}, \dots, \tau_r) = G_r(z+\tau_k, \tau_0, \dots, \tau_r)^{-1}.\end{equation} 
Starting from any $\tau_0, \dots, \tau_r \in \cc - \rr$, we may define $J = \{ 0 \leq j \leq r \setseparator \Im(\tau_j) < 0\}$. Then using the inversion property \refp{inversionpropertyalgo} for each $k \in J$ gives:
$$G_r(z, \tau_0, \dots, \tau_r) = G_r\left(z - \sum_{k \in J} \tau_k, \tau_0', \dots, \tau_r'\right)^{(-1)^{\cardinalshort{J}}} $$
where $\tau'_j = \tau_j$ if $j \not\in J$ and $\tau'_j = -\tau_j$ if $j \in J$. In any case, this gives $\Im(\tau'_j) > 0$. The number of reorientation steps to perform is at most $r+1$.

The translation step consists of bringing the elliptic variable $z$ into the center strip defined by $0 < \Im(z) < \sum_{j = 0}^r \Im(\tau_j)$. For this step we use recursively the pseudo-periodicity property:
\begin{equation}\label{pseudoperiodicityalgo} 
G_r(z+\tau_k, \tau_0, \dots, \tau_r) = G_{r-1}(z, \tau_0, \dots, \omitvar{\tau_k}, \dots, \tau_r) \cdot G_r(z, \tau_0, \dots, \tau_r).
\end{equation} 
Doing so will require the computation of lower degree functions, that is, until we reach the case $r = 0$ in which case the computation of $G_0 = \theta$ can be done quickly using \refp{jacobitriple}. In the general case, the number of translation steps required to reach the center strip for all involved functions, including the new ones appearing at each step is $O(\prod_{k = 1}^r (|\Im(z)|/\Im(\tau_k)))$. The proof of this statement may be done by induction, but we shall skip it as we now focus on the specific computation of the conjectural higher elliptic unit which doesn't require many translation steps. Indeed, we may actually prove that each of the terms in the product \refp{defukb} requires either $0$ or $r = n-2$ translation steps to be computed:

\begin{general}{Proposition}
Fix a permutation $\rho \in \goth{S}_{n-2}$. Set $J_{\rho} = \{ 1 \leq j \leq n-1 \setseparator \Im(\tau_{j, \rho}) < 0\}$. The term 
$$ T_{\rho} = \frac{G_{n-2}\left(z_{\rho}, \tau_{1, \rho}, \dots, \tau_{n-1, \rho}\right)^N}{G_{n-2}\left(Nz_{\rho}, N\tau_{1, \rho}, \dots, N\tau_{n-1, \rho}\right)}$$
with the parameters $z_{\rho}, \tau_{1, \rho}, \dots, \tau_{n-1, \rho}$ defined in \refp{defzparam} and \refp{deftauparam} is equal to 
$$T_{\rho} = \left(\frac{G_{n-2}\left(z'_{\rho}, \tau'_{1, \rho}, \dots, \tau'_{n-1, \rho}\right)^N}{G_{n-2}\left(Nz'_{\rho}, N\tau'_{1, \rho}, \dots, N\tau'_{n-1, \rho}\right)}\right)^{(-1)^{\cardinalshort{J_{\rho}}}}$$
with 
\begin{align*}
z'_{\rho} &= z_{\rho} -\sum_{j \in J_{\rho}} \tau_{j, \rho}, \\ 
\tau'_{j, \rho} &= \begin{cases} \tau_{j, \rho} & \text{ if } j\not\in J_{\rho} \\ -\tau_{j, \rho} & \text{ if } j \in J_{\rho} \end{cases} \\
0 & \leq \Im(z'_{\rho}) \leq \sum_{j = 1}^{n-1} \Im(\tau'_{j, \rho}). 
\end{align*}
As a consequence, if $0 < \Im(z'_{\rho}) < \sum_{j = 1}^{n-1} \Im(\tau'_{j, \rho})$ then $T_{\rho}$ requires $0$ translations steps to be computed. Otherwise, $T_{\rho}$ requires $n-2$ translation steps to be computed.
\end{general}

\begin{proof}
It follows from the definition of the evaluation parameters $\tau_{1, \rho}, \dots, \tau_{n-1, \rho}$ that the family $1, \tau_{1, \rho}, \dots, \tau_{n-1, \rho}$ is a $\qq$-basis of $\kk$. By assumption, $\kk$ contains no subfield, which implies that $\kk \cap \rr = \qq$. Thus the parameters $\tau_{1, \rho}, \dots, \tau_{n-1, \rho}$ lie in $\cc - \rr$. The formula for $T_{\rho}$ in terms of the $z'_{\rho}, \tau'_{1, \rho}, \dots, \tau'_{n-1, \rho}$ is obtained via the reorientation step. The last claim follows from the fact that $z_{\rho} \in \qq$ so 
$$0 \leq \Im(z'_{\rho}) = \sum_{j \in J_{\rho}} \Im(\tau'_{j, \rho}) \leq \sum_{j = 1}^{n-1} \Im(\tau'_{j, \rho}). $$
In particular, if $\emptyset \neq J_{\rho} \subsetneq \{1, \dots, n-1\}$ then $z'_{\rho}$ is already in the center strip. Otherwise, using $n-2$ successive translation steps all the way down to $G_0$ one has the expression
$$T_{\rho} = \prod_{r = 0}^{n-2}\left(\frac{G_{n-2-r}\left(z'_{\rho} + \tau'_{r+1, \rho}, \tau'_{r+1, \rho}, \dots, \tau'_{n-1, \rho}\right)^N}{G_{n-2-r}\left(N(z'_{\rho}+\tau'_{r+1, \rho}), N\tau'_{r+1, \rho}, \dots, N\tau'_{n-1, \rho}\right)}\right)^{(-1)^{n-2+r+\cardinalshort{J_{\rho}}}}$$
when $\Im(z'_{\rho}) = 0$ and 
$$T_{\rho} = \prod_{r=0}^{n-2}\left(\frac{G_{n-2-r}\left(z'_{\rho} - \sum_{j = 1}^{r+1}\tau'_{j, \rho}, \tau'_{r+1, \rho}, \dots, \tau'_{n-1, \rho}\right)^N}{G_{n-2-r}\left(N(z'_{\rho} - \sum_{j = 1}^{r+1} \tau'_{j, \rho}), N\tau'_{r+1, \rho}, \dots, N\tau'_{n-1, \rho}\right)}\right)^{(-1)^{\cardinalshort{J_{\rho}}}}  $$
when $\Im(z'_{\rho}) = \sum_{j = 1}^{n-1} \Im(\tau'_{j, \rho})$. Each of these expressions contains only terms with the elliptic variable in the center strip.
\end{proof}

It follows from this Proposition that the computation time to check the conjecture is entirely dedicated to the computation of sums of the shape \refp{sumcomputations}. Let us say a few words on the number of terms we must compute in this sum. Let $y_{\rho} = 2\pi.\min(\Im(z'_{\rho}), \Im(\tau'_{1, \rho}) + \dots + \Im(\tau'_{n-1, \rho}) - \Im(z'_{\rho}))$ be the distance from $z'_{\rho}$ to the boundary of the center strip. Then, to compute one of the terms $G_{n-2}(z'_{\rho}, \tau'_{1, \rho}, \dots, \tau'_{n-1, \rho})$ appearing in the expression for $u_{k, \goth{b}}$ with precision $\delta > 0$ we need to compute $O(\log(\delta)/y_{\rho})$ terms in at most $(n-2)$ sums of the form \refp{sumcomputations} (the computation of $G_0 = \theta$ is negligible). The main difficulty in performing the computations comes from the fact that, in practice, the imaginary parts of the parameters $\tau_{1, \rho}, \dots, \tau_{n-1, \rho}$ are very small, making the parameter $y_{\rho}$ small as well.

\subsection{Sieving for optimal settings}

Let us now describe how to obtain optimal examples for our generalised elliptic units. First, we explain how to obtain number fields for which the optimal conditions are satisfied for some suitable choice of class field modulus $\goth{f}$ and for some choice of fundamental units for $\opcf$. Then, we explain how to compute the higher elliptic units for such a number field.

To search for optimal examples, we shall assume that we are provided with some iterator over integral polynomials of fixed degree $n$. This iterator might be very broad, describing all polynomials of degree $n$ with coefficients between $-A$ and $A$ for some $A>0$ or describe particular polynomials such as the pure cubic polynomials $x^3-m$ for $1 \leq m \leq A$. To gain time, it might be useful to use an iterator built from a list of all pre-computed polynomials defining a number field of degree $n$ with exactly one complex place with discriminant less than $A$ for some $A> 0$, with possibly additional conditions on the class number for instance. Such lists may be found on the LMFDB \cite{lmfdb} for small degrees. \bigskip

\noindent \textbf{Algorithm 1:} Search for optimal examples.
\newline \noindent \textbf{Input:} a degree $n \geq 3$ and an iterator $g : \bb{N} \to \zz_n[X] = \{ P \in \zz[X]\setseparator \mathrm{deg}(P) = n\}$.
\newline \noindent \textbf{Output:} a list $L$ of vectors $[\kk, \goth{f}_{list}, \eps_{list}]$ where $\kk$ is a number field of degree $n$ with exactly one complex place, $\goth{f}_{list} = [\goth{f}_1, \dots, \goth{f}_l]$ contains class field moduli and $\eps_{list} = [\eps_1, \dots, \eps_{n-2}]$ contains fundamental units for $\opck$ such that the data $(\kk, \goth{f}_i, \eps_{list})$ satisfies optimal conditions for the construction of higher elliptic units.
\newline \noindent Initialise a list $L \leftarrow List()$;
\newline \noindent Then, for each polynomial $P$ described by the iterator $g$ do the following:
\begin{enumerate}
\item If $P$ is reducible, move on to the next polynomial. 
\item Compute the number $r$ of real roots of $P$. If $r \neq n-2$, move on to the next polynomial.
\item Let $\kk = \qq[X]/(P)$. Compute a $\zz$-basis of the maximal order $\ok$ of $\kk$. 
\item Compute a set $\eps_1, \dots, \eps_{n-2}$ of LLL-reduced fundamental units for $\opck$ (see for instance the \textbf{bnfinit} command in Pari/GP).
\item Compute the contents $\tlambda_{\rho}$ and the linear forms $\tilde{a}_{\rho}$ for all permutation $\rho \in \goth{S}_{n-2}$. If at any point a content $\tlambda_{\rho} > 1$ is found, move on to the next polynomial.
\item Compute the HNF representation of the associated overflow ideals $\tildeD_{\rho}$ using a variant of [\hspace{1sp}\cite{cohen}, Proposition 4.8.19] (see \cite{thesis} for more details). 
\item Compute the prime factorisation of the ideal $J(\kk) = \sum_{j = 1}^{n-2}(\eps_j -1)\ok$ (see \refp{defjk}) and keep only ideals $\goth{f}_1, \dots, \goth{f}_l$ dividing $J(\kk)$ which are coprime to $\prod_{\rho}\tildeD_{\rho}$ and for which $\ok/\goth{f}_i$ is cyclic. To avoid examples where the higher elliptic unit is expected to be a root of unity, one may decide to keep only those modulus for which $\kk^{+}(\goth{f}_i)$ is totally complex. Note that the signature of $\kk^{+}(\goth{f}_i)$ can be read off from the ray class group $\classgroup{\goth{f}_i}$ (see the \textbf{bnrdisc} command in Pari/GP). If no such prime factor may be found, move on to the next polynomial. Otherwise, add the data $[\kk, [\goth{f}_1, \dots, \goth{f}_l], [\eps_1, \dots, \eps_{n-2}]]$ to the list $L$.
\end{enumerate}

Let us make a few remarks on Algorithm 1. In practice, as the optimal conditions depend on the choice of fundamental units $\eps_1, \dots, \eps_{n-2}$ for $\opck$, it may sometimes be worth it to try more than one choice for $\eps_1, \dots, \eps_{n-2}$ by performing small base changes around an LLL-reduced choice of fundamental units. As a second remark, the optimal conditions we described (see section \ref{sectionhypotheses}) do not guarantee that the computations will be feasible in a reasonable amount of time. Indeed, the computation time for the higher elliptic units $u_{k, \goth{b}}$ depends on the position of the evaluation parameters $\tau_{1, \rho}, \dots, \tau_{n-1, \rho}$. In particular, these evaluations parameters have denonimators $\ell = q\cdot N \cdot \ttt_{\rho} \cdot p$ where the helper ideal $\goth{H}$ as norm $p^{n-1}$. Thus, in general, high values of $\ttt_{\rho}$ may lead to unfeasible computations. For instance, the pure cubic field $\qq(z = e^{2i\pi/3}\sqrt[3]{93})$ has fundamental unit $\eps = 15001z^2 - 64428z - 16022$ and the corresponding value for $\ttt$ is $648833101994018933678601952991$ while $\tlambda = 1$. In practice, we modify Algorithm 1 to discard potential examples with high values of $\ttt_{\rho}$.

Using this algorithm we have already obtained more than $10000$ complex cubic fields, several thousands of quartic and quintic fields and a few hundreds of degree $6$ fields satisfying these optimal conditions. The best scenario we found for degree $7$ fields was a case where $116$ out of the $120$ conditions $\tlambda_{\rho} = 1$ for $\rho \in \goth{S}_{5}$ were satisfied, and in that case the values of $\ttt_{\rho}$ were too large to hope we could compute the higher elliptic units.

\subsection{Computing higher elliptic units in optimal settings}

We now move on to the computation of the conjectural elliptic units, provided that we have selected an optimal setting beforehand. Let us indeed suppose that we are given a number field $\kk$ of degree $n \geq 3$ with exactly one complex place, a degree one prime $\goth{f}$ above the prime $q$ in $\kk$ and a set of fundamental units $\eps_1, \dots, \eps_{n-2}$ for $\opck = \opcf$. The following algorithm describes the computation of higher elliptic units for this optimal setting: \bigskip

\noindent \textbf{Algorithm 2:} Computation of higher elliptic units in optimal cases .
\newline \noindent \textbf{Input:} $\kk$ a number field of degree $n \geq 3$ with exactly one complex place, $\goth{f}$ a degree one prime in $\kk$ such that $\classfieldf$ is totally complex and a set of fundamental units $\eps_1, \dots, \eps_{n-2}$ for $\opck = \opcf$ as given by Algorithm 1.
\newline \noindent \textbf{Output:} a list of conjectural elliptic units $u_{k, \goth{b}}$ given by \refp{defukb}.
\begin{enumerate}
\item Compute the overflow ideals $\tildeD_{\rho}$ following Algorithm 1. Set $\tildeD = \prod_{\rho} \tildeD_{\rho}$.
\item Compute the ray class group $\classgroupf$. 
\item Find a prime ideal $\goth{a}$ of prime norm $N$ such that $N$ is coprime to $\goth{f} \cdot \tildeD$. 
\item Fix a finite set $R$ of integral ideals coprime to $N \cdot \goth{f}\cdot \tildeD$ to represent the narrow Hilbert class group of $\kk$.
\item For each $\goth{b} \in R$, compute a helper ideal $\goth{H}$ and the base points $h_{\rho}$ as given by Definition \ref{defcompatible}. Compute the evaluation parameters $z_{\rho}, \tau_{1, \rho}, \dots, \tau_{n-1,\rho}$ as described in \refp{defzparam} and \refp{deftauparam}.
\item For each $\goth{b} \in R$ and each $k \in (\zsz{q})^{\times}$, compute $u_{k, \goth{b}}$ for all possible choices of signs $\epsilon_{\rho}, \epsilon'_{\rho}$ and identify the correct one using the expected Kronecker limit formula \refp{klfconjecture}. 
\end{enumerate}

\textbf{Remark:} In practice, to ``find an ideal'', we iterate over all integral ideals of $\kk$ by increasing norm until the set of representatives $R$ is complete and until we find a suitable smoothing ideal $\goth{a}$. Steps 3 and 4 may be interverted, however, we usually fix the smoothing ideal first as the computation time depends linearly on $N$ whereas it is independent of the choice of ideals $\goth{b}$. Fixing $N$ first usually allows us to keep $N$ relatively small.

\section{Numerical evidence to support the conjecture}\label{sectionnumerical}

In this section, we provide numerical examples to support our conjecture\footnotemark[2]\footnotetext[2]{The computations may be checked using the files provided at:\newline https://plmlab.math.cnrs.fr/pmorain/computations\_of\_higher\_elliptic\_units} in optimal settings. They may be computed with high precision in a short amount of time. In what follows, we will give computation times for $1000$ digits precision on a personal computer. Computations were performed using number fields found in the LMFDB database \cite{lmfdb} and carried out using the computer algebra system PARI/GP \cite{parigp}, making extensive use of algebraic number theory tools it provides. To understand the result of the computations, we use methods relying on the LLL algorithm to search for linear relations among values. Namely, we use the commands \textbf{lindep} to test if the value $\log|u_{k, \goth{b}}|^2$ obtained from our computions is close to the expected value for $N\zeta'_{\goth{f}}([k\goth{b}], 0) - \zeta'_{\goth{f}}([k\goth{a}\goth{b}], 0)$ and $\textbf{algdep}$ to test if the value obtained for $u_{k, \goth{b}}$ is close to some algebraic integer in an abelian extension of $\kk$. Another possibility to interpret the results is to directly compute the polynomial
$$\prod_{k \in \zsz{q}^{\times}}\prod_{\goth{b} \in \classgroupplusk} (X - u_{k, \goth{b}})$$
which should lie in $\ok[X]$ and identify its coefficients as elements of $\kk$ using $\textbf{lindep}$. In any case, we may check afterwards that any obtained candidate polynomial defines the number field $\classfieldf$ above $\qq$ or $\kk$.

In what follows, we will define our fields as $\kk = \qq(z)$ where $z$ is the unique complex root of some polynomial $P \in \qq[X]$ lying in the upper half-plane. We provide examples of computations of higher elliptic units above number fields of degree 3, 4, 5 and 6 with exactly one complex place (in optimal settings). In all cases, we checked that the values we obtain correspond to the embeddings of a unit inside $\kk^{+}(\goth{f})$ satisfying the expected Kronecker limit formula \refp{klfconjecture}. Examples in non-optimal settings are presented in \cite{thesis}.

\stepcounter{example}
\subsection{Example \exc: a detailed cubic example with $q = 11$}

We start this section of examples with a detailed cubic example with trivial narrow Hilbert class group and $q = 11$. Let us now detail Algorithm 1 for the integral polynomial $P = x^3 - x^2 + 5x - 2$. It is irreducible and its discriminant is $D = -411$ so it has one real root $x_0$ and two complex roots $\zexc, \overline{\zexc}$ such that $\Im(z) > 0$. We then put $\kk = \qq(\zexc)$ and compute a $\zz$-basis $B = [1, z, z^2]$ of $\ok$. A fundamental unit for $\opck$ is $\eps = \zexc^2+2\zexc-1$. The linear form $\tilde{a}$ is defined by $\tilde{a}(\beta_0 + \beta_1\zexc + \beta_2\zexc^2) = 2\beta_2-\beta_1$ and the content $\tlambda$ is equal to $1$. The overflow ideal $\tildeD$ associated to $\tilde{a}$ is the unique ideal of norm $31$ in $\ok$ and it is generated by $z-3$. We then compute the ideal $J(\kk) = (\eps-1)\ok = \goth{P}_2\goth{P}_3\goth{P}_{11}$ (see \refp{defjk}) where $\goth{P}_2$ is the unique ideal of norm $2$, $\goth{P}_{11}$ is the unique ideal of norm $11$ and $3\ok = \goth{P}_3\goth{P}^{\prime2}_3$ in $\ok$. Any divisor of $J(\kk)$ which is neither $\ok$ nor $\goth{P}_2$ gives a totally complex class field. For the purposes of this example, we choose $\goth{f}$ to be the ideal $\goth{P}_{11}$ of norm $q = 11$. 

We shall now describe Algorithm 2 applied to this example. The overflow ideal has already been computed: $\tildeD$ is the ideal of norm $31$ in $\ok$. The narrow ray class group at $\goth{f}$ is given by $\classgroupf \simeq (\zsz{11})^{\times} \simeq \zsz{10}$. Thus, we shall fix $\goth{b} = \ok$ and vary $k$ in $(\zsz{11})^{\times}$. The smoothing ideal $\goth{a}$ is chosen to be the unique ideal of norm $17$ in $ \ok$ (we checked that $17$ is indeed coprime to $q \cdot \tildeD = 11 \cdot 31$. Following \refp{defcompatible} we may choose the base point $h = -121z^2 + 44z - 407$ which is a generator of the ideal $17\cdot11\cdot\goth{D}\cdot(\goth{a}\goth{b})^{-1}$, that is, $m = 1$ and $\goth{H} = \ok$. Our construction then gives the evaluation parameters:
$$\tau = 12\zexc^2 - 7\zexc + 1310, ~~ \sigma = \zexc^2 + 2\zexc + 383 $$
so that $u_{k, \ok}$ is given by the $17$-smoothed evaluation:
$$u_{k, \ok} =  \frac{\Gamma\left(\frac{k}{11}, \frac{\tau}{5797}, \frac{\sigma}{5797}\right)^{-17}}{\Gamma\left(\frac{17k}{11}, \frac{17\tau}{5797}, \frac{17\sigma}{5797}\right)^{-1}}\hfill\\$$
for any $1 \leq k \leq 10$. These complex numbers are computed with high precision to be close to the ten roots of the following relative palindromic polynomial:
$$P = x^{10} + c_1(z)x^9 + c_2(z) x^8 + c_3(z)x^7 + c_4(z) x^6 + c_5(z) x^5 + c_4(z) x^4 + c_3(z)x^3 + c_2(z) x^2 + c_1(z) x + 1 $$
defining the extension $ \kk^{+}(\goth{f})/ \kk$, where
\begin{align*}
c_1(z) & = -882z^2-310466z+130708 \\
c_2(z) & = 23602645z^2 - 2290686z - 3210066 \\
c_3(z) & = -3650154729z^2 + 5329546969z - 1595637486\\
c_4(z) & = -67756252890z^2 + 210542793317z - 76551497367\\
c_5(z) & = -96836154271z^2 + 989098538346z - 398787019925.
\end{align*}
One can also check the Kronecker limit formula of Bergeron, Charollois and Garc\'ia \cite{BCG} as for instance:
$$17\cdot \zeta'_{\goth{f}}([\ok], 0) - \zeta'_{\goth{f}}([\goth{a}], 0) \approx \log|u_{1, \ok}|^2 \approx 3.3591950...$$

\stepcounter{example}
\subsection{Example \exc: a cubic example with a large class group}

Let $\zexc$ be the complex root of the polynomial $x^3-65$ in the upper half-plane. The number field $ \kk = \qq(\zexc)$ has class number $18$. We set $\goth{f}$ to be the unique prime ideal of norm $q = 3$; this ideal satisfies $\goth{f}^3 = (3)$ and $(\kk, \goth{f})$ is an optimal setting. The corresponding narrow ray class group is $\classgroupf \simeq \zsz{6}\times \zsz{6}$. The group $\opcf$ is generated by $\eps = z-4$ and the overflow ideal $\tildeD$ is trivial. We give the computations for $\goth{b} = \ok$, $\goth{b} = \goth{P}_{59}$ where $\goth{P}_{59}$ is the unique integral ideal of norm $59$ in $ \kk$ and $k = 1, 2$, that is 4 out of the 36 class in $\classgroupf$. The smoothing ideal $\goth{a}$ may be chosen as the unique prime ideal of norm $N = 5$ in $ \ok$. Our method gives:
\begin{align*}
u_{k, \ok} & = \frac{\Gamma\left(\frac{k}{3}, \frac{-\zexc^2 - 4\zexc}{195}, \frac{-\zexc+65}{195}\right)^{5}}{\Gamma\left(\frac{5k}{3}, \frac{5(-\zexc^2 - 4\zexc)}{195}, \frac{5(-\zexc+65)}{195}\right)}  \approx \begin{cases} -1.6691052... +i\cdot5.7493283... \text{ for } k = 1\\ -0.0465701...-i\cdot 0.1604134... \text{ for } k = 2\end{cases}\hfill\\
u_{k, \goth{P}_{59}} & = \frac{\Gamma\left(\frac{k}{3}, \frac{\zexc^2 + 4\zexc + 1050}{615}, \frac{\zexc-95}{615}\right)^{-5}}{\Gamma\left(\frac{5k}{3}, \frac{5(\zexc^2 + 4\zexc + 1050)}{615}, \frac{5(\zexc-95)}{615}\right)^{-1}}  \approx \begin{cases}  0.0344135... - i\cdot0.0123218... \text{ for } k = 1\\ 25.7563258...+i\cdot 9.2221524... \text{ for } k = 2 \end{cases}\hfill
\end{align*}
where $195 = 3\cdot5\cdot13$, $615 = 3\cdot5\cdot41$, the factors $13$ and $41$ arising from the helper ideals above $13$ and $41$ respectively. We may compute the remaining 32 out of 36 values $u_{k, \goth{b}}$ attached to the 36 classes in $\classgroupf$ and identify the palindromic polynomial:
$$\prod_{k, \goth{b}} (X - u_{k, \goth{b}}) = X^{36}-(a\zexc^2 +b\zexc +c)X^{35} + \dots + 1 \in \mathcal{O}_{ \kk}[X]$$
where for instance $a = -5967373310133, b = 769211619985, c= 93377174024326$. This polynomial defines a relative equation of the class field $ \kk^{+}(\goth{f})$ above $ \kk$ and we identify the rest of its coefficients in $\mathcal{O}_{ \kk}$ using the \textbf{lindep} command. The computation time for $1000$ digits and for all of the 36 computations is 8 minutes, which gives 13 seconds per individual computation on average. For this specific computation, one may remark that for all classes in $\classgroupf$, $u_{2, \goth{b}} = u_{1, \goth{b}}^{-1}$ which is quite typical in the cubic case when $q = 3$.

It is interesting to note that this field belongs to a special family parametrized by $P_m = x^3 - 2^{3m} - (-1)^m$ for $m \geq 0$ which behaves nicely compared to other pure cubic fields for the specific modulus $\goth{f}$ satisfying $\goth{f}^3 = (3)$. Putting $\bb{L}_m = \qq(z_m)$ where $z_m = e^{2i\pi/3} \sqrt[3]{2^{3m} + (-1)^m}$, it seems that the positive fundamental unit in $\bb{L}_m$ is given by $\eps_m = (-1)^m(z_m-2^m)$ and its inverse by $\eps_m^{-1} = z_m^2 + 2^mz_m + 2^{2m}$. When the discriminant $D_m = 2^{3m} + (-1)^m$ of $x^3 - 2^{3m} - (-1)^m$ is cube-free (at least for $0 \leq m \leq 100$ except for $m = 49, 50$), because $D_m \not\equiv \pm 1 \mod 9$, the field $\bb{L}_m$ is a so-called pure cubic field of the first kind and an integral basis of $\mathcal{O}_{\bb{L}_m}$ is given by $(1, z_m, z_m^2/s_m)$ where $D_m = r_m s_m^2$ with $r_m, s_m$ square-free and relatively coprime.  When $D_m$ is also square-free (at least for $0 \leq m \leq 100$ except for $m = 7, 10, 21, 26, 30, 35, 63, 68, 70, 77, 78, 90, 91$) we get optimal conditions to perform the computations as the overflow ideal governing our computations is trivial. However, the class group grows very rapidly in this family (for instance: $h(\bb{L}_{11}) = 2^2\cdot3^5\cdot5\cdot 3191$ and
$h(\bb{L}_{20}) = 2^6\cdot3^4\cdot5\cdot11\cdot19\cdot79\cdot 863\cdot18047$) so the results of our computations become more difficult to interpret as more precision is required to find linear dependence relations. The following table gives information about the first few of these fields:

\begin{table}[ht!]
\begin{center}
\noindent \begin{tabular}{| C | C | C | C | C | C | C |}
\hline
m & 0 & 1 & 2 & 3 & 4 & 5\\\hline
&&&&&&\\[-1em]
P_m & x^3-2 & x^3-7 & x^3-65 & x^3 -511 & x^3 - 4097 & x^3-32767\\\hline
&&&&&&\\[-1em]
h(\bb{L}_m) & 1 & 3 & 2\cdot3^2 & 2^2\cdot3^3 & 2^3\cdot3^4 & 2^3\cdot3^3\cdot5^2 \\
\hline
\end{tabular} 
\label{tablecubfamily}
\end{center}
\end{table}

We have computed the higher elliptic units for the fields $\bb{L}_m$, $m = 0, 1, 2, 3, 4$, but already for $m = 5$ the relative polynomial $\prod_{k, \goth{b}}(X - u_{k,\goth{b}})$ has expected degree $10800$ and huge coefficients, thus we would need a lot of precision in the computation of the higher elliptic units $u_{k, \goth{b}}$ to correctly identify it.

Lastly, it would be interesting to understand if for $m \geq 100$ there are infinitely many fields in this family with squarefree $D_m$ as our construction would then conjecturally allow to construct algebraic units of very high degree.

\stepcounter{example}
\subsection{Example \exc: one of the simplest quartic examples}\label{sectionsimplestquartic}

Let us now illustrate the conjecture with one of the simplest examples where the base field $\kk$ is a quartic number field. Let $\zexc$ be the complex root of the polynomial $x^4 -6x^3-x^2-3x+1$ lying in the upper half-plane. The number field $\kk = \qq(\zexc)$ has class number $1$ and exactly one complex place. We choose $\goth{f}$ to be the unique ideal of norm $q = 2$ in $\ok$. We check that $\opck = \opcf$ and we fix a set $\eps_1, \eps_2$ of fundamental units for $\opck$ given by:
$$ \eps_1 = \frac{-2\zexc^3 + 13\zexc^2 - \zexc + 3}{7},~\eps_2 = \frac{-5\zexc^3 + 29\zexc^2 + 15\zexc + 18}{7}.$$
The associated overflow ideals $\tildeD_1, \tildeD_2$ are trivial. The narrow ray class group $\classgroupf \simeq \zsz{2}$ is isomorphic to the narrow Hilbert class group and it is represented by the ideals $\ok$ and $\goth{P}_3$ where $\goth{P}_3$ is the unique prime ideal of norm $3$. We may fix the smoothing ideal $\goth{a}$ to be the unique prime ideal of norm $N = 13$ in $\ok$. Our method gives the evaluation parameters
\begin{align*}
&\tau = 5z^3 - 29z^2 - 15z + 87, && \tau' = 2z^3 - 13z^2 + z - 24 = - \rho \\
&\sigma = 6z^3 - 39z^2 + 10z + 47, && \sigma' = 5z^3 - 29z^2 - 15z + 87 = \tau  \\
&\rho = -2z^3 + 13z^2 - z + 24, && \rho' = -2z^3 + 13z^2 + 6z + 143
\end{align*}
so that the complex numbers $u_{1, \ok}$ and $u_{1, \goth{P}_3}$ are given by:
$$u_{1, \ok} = \frac{G_2\left(\frac{1}{2}, \frac{\tau}{182}, \frac{\sigma}{182}, \frac{\rho}{182}\right)^{13}}{G_2\left(\frac{13}{2}, \frac{13\tau}{182}, \frac{13\sigma}{182}, \frac{13\rho}{182}\right)} \times \frac{G_2\left(\frac{1}{2}, \frac{\tau'}{182}, \frac{\sigma'}{182}, \frac{\rho'}{182}\right)^{13}}{G_2\left(\frac{13}{2}, \frac{13\tau'}{182}, \frac{13\sigma'}{182}, \frac{13\rho'}{182}\right)} = u_{1, \goth{P}_3}^{-1}.$$
The complex number $u_{1, \ok}$ is computed to 1000 digits of precision and is found to be close to the root $ \approx 4.1210208... - i\cdot5.0617720...$ of the polynomial 
$$x^8 - 7x^7 + 33x^6 + 49x^5 + 17x^4 + 49x^3 + 33x^2 - 7x + 1$$
which defines $\kk^{+}(\goth{f})$ over $\qq$. The computation time for $1000$ digits is 6 seconds. We may also check the conjectured Kronecker limit formula \refp{klfconjecture} up to 1000 digits as:
$$13\cdot \zeta'_{\goth{f}}([\ok], 0) - \zeta'_{\goth{f}}([\goth{a}], 0) \approx \log\left|u_{1, \ok}\right|^2 \approx 3.7519563...$$

\stepcounter{example}
\subsection{Example \exc: a quartic example with $q = 7$}\label{quarticq7}

We now discuss a quartic example with $q = 7$ and non-trivial class group. Let $\zexc$ be the complex root of the polynomial $x^4 - x^3 - 4x^2 - 11x + 16$ lying in the upper half-plane. Then $ \kk = \qq(\zexc)$ has class number $2$ and narrow Hilbert class group $\classgroupplus{ \kk} \simeq (\zsz{2})^2$. Let us fix the modulus $\goth{f}$ to be the unique prime ideal of norm $q = 7$ in $ \ok$. The corresponding narrow ray class group is $\classgroupf \simeq \zsz{12} \times \zsz{2}$. If we choose the fundamental units
$$ \eps_1 = \zexc -1,~\eps_2 = -\zexc + 3 $$
for $\opc{\goth{f}} = \opck$, the overflow ideals governing our computation are $\tildeD_1 = \tildeD_2 = \goth{P}_5^2$ where the factorisation of $(5)$ in terms of prime ideals is given by $(5) = \goth{P}_5^3 \goth{P}_5'$. Let us fix the smoothing ideal $\goth{a} = 19\ok + (\zexc + 7)\ok$ to be one of the two ideals of norm $N = 19$ in $\ok$. For the ideal $\goth{b} = \ok$, our method gives the evaluation parameters
\begin{align*}
&\tau = \zexc^3 - 4\zexc + 4172, && \tau' = \zexc^3 + 2\zexc^2 + 2\zexc - 2344 \\
&\sigma = \zexc^2 + 3\zexc - 3258, && \sigma' = -\zexc^2 - 3\zexc + 3258 = -\sigma \\
&\rho = -\zexc - 1356, && \rho' = -\zexc - 1356 = \rho
\end{align*}
which we use to evaluate $u_{k, \ok}$. Indeed, we may compute for $k = 1, 2, 3, 4, 5, 6$:
\begin{align*}
u_{k, \ok}&= \frac{G_2\left(\frac{k}{7}, \frac{\tau}{1330}, \frac{\sigma}{1330}, \frac{\rho}{1330}\right)^{19}}{G_2\left(\frac{19k}{7}, \frac{19\tau}{1330}, \frac{19\sigma}{1330}, \frac{19\rho}{1330}\right)} \times \frac{G_2\left(\frac{k}{7}, \frac{\tau'}{1330}, \frac{\sigma'}{1330}, \frac{\rho'}{1330}\right)^{19}}{G_2\left(\frac{19k}{7}, \frac{19\tau'}{1330}, \frac{19\sigma'}{1330}, \frac{19\rho'}{1330}\right)}\\
u_{k, \ok}& \approx \begin{cases} -10781954.9966390...-i\cdot 1533647.9705555... \text{ for } k = 1, 6 \\
5.9031254...-i\cdot 4.7220508... \text{ for } k = 2, 5 \\
0.1823482...-i\cdot 0.1028436... \text{ for } k = 3, 4 \end{cases}
\end{align*}
to high precision, where $1330 = 2\cdot5\cdot7\cdot19$. The fact that some of the values are equal shows that the unit we compute actually belongs to a subextension $\bb{L}/ \kk$ of $ \kk^{+}(\goth{f})/ \kk$ of index $2$. We may compute similarly the remaining $18$ values corresponding to $k = 1, 2, 3, 4, 5, 6$ and to the three non-trivial classes in $\classgroupplus{ \kk}$, and check that $u_{6-k, \goth{b}} = u_{k, \goth{b}}$ in all cases. We then identify the relative palindromic polynomial
$$P = \prod_{k = 1}^3 \prod_{\goth{b}} (X - u_{k, \goth{b}}) = X^{12} + (a\zexc^3 + b\zexc^2 + c \zexc + d \zexc)X^{11} + \dots \in \ok[X] $$
where for instance
\begin{align*}
a & = 2336913308/5\\
b & = 6736333629/5\\
c &= -57839312901/5\\
d &= 50793489608/5
\end{align*}
which is coherent with the fact that a $\zz$-basis of $\ok$ is given by $B = [1, \zexc, (\zexc^3 - 2\zexc^2 - 2\zexc - 4)/5, (2\zexc^3 + \zexc^2 - 9\zexc - 23)/5]$. We check that $P$ defines an extension $\bb{L}$ of $\kk$ of degree $48$ over $\qq$ which is a subextension of $\classfieldf/\kk$. From the point of view of Stark's conjecture, since the field $ \classfieldf$ contains the roots of unity of order $42$, the conjectural units we compute are essentially $12$th powers of the conjectural Stark units and we may construct the total class field $\kk^{+}(\goth{f})$ by extracting roots in a suitable manner. The main difficulty in treating this kind of examples with a large ray class group is that a lot of precision is required to analyse the result and identify the coefficients of the relative polynomial at the end, since it requires finding linear dependence relations in vector spaces of high dimension.

\stepcounter{example}
\subsection{Example \exc: a quintic example}
We now discuss one of our simplest quintic examples, which is already quite long to write down as it requires to describe a product of $(5-2)! = 6$ smoothed $G_3$ functions, each involving $4$ parameters (see Table \ref{tablequintic}). Let $\zexc$ be the complex root of the polynomial $x^5 - x^4 - x^3 - 2x^2 + x + 1$ lying in the upper half-plane. The number field $ \kk = \qq(\zexc)$ has narrow class number $1$. We choose $\goth{f}$ to be the unique prime ideal of norm $q = 3$. The corresponding narrow ray class group is $\classgroupf \simeq \zsz{3}^{\times} \simeq \zsz{2}$. If we fix the fundamental units
$$ \eps_1 = \zexc^4 - 2\zexc^3 - \zexc + 3,~\eps_2 = 2\zexc^4 - 2\zexc^3 - 3\zexc + 3,~\eps_3 = 2\zexc^4 - 3\zexc^3 - 4\zexc + 4 $$
for $\opc{\goth{f}} = \opck$ and the order $\{ \mathrm{Id}, (32), (21), (231), (312), (31)\}$ of $\goth{S}_3$, then the overflow ideals governing our construction are integral ideals of squarefree norm: 
\begin{align*}
\norm{\goth{D}_1} & = 7\cdot 37 \cdot 137  && \norm{\goth{D}_2} = 31\cdot53 && \norm{\goth{D}_3} = 491\\
\norm{\goth{D}_4} &= 107  && \norm{\goth{D}_5} = 1 && \norm{\goth{D}_6} = 145637.
\end{align*}
We may choose $\goth{b} = (1)$, $k = 1,2$ and $\goth{a}$ the unique prime ideal of norm $N = 11$ in $ \ok$. We then compute six quotients: 
$$v_j = \frac{G_3\left(\frac{1}{3}, \frac{\tau_j}{l_j}, \frac{\sigma_j}{l_j}, \frac{\rho_j}{l_j}, \frac{\varpi_j}{l_j}\right)^{11}}{G_3\left(\frac{11}{3}, \frac{11\tau_j}{l_j}, \frac{11\sigma_j}{l_j}, \frac{11\rho_j}{l_j}, \frac{11\varpi_j}{l_j}\right)} $$
where the parameters $\tau_j, \sigma_j, \rho_j, \varpi_j$ are given in Table \ref{tablequintic} below and we have defined for ease of presentation the levels $l_j = 3 \cdot 11 \cdot \norm{\goth{D}_j}$, that is:
$$l_1 =1170939, ~l_2 = 54219, ~l_3 = 16203, ~l_4 = 3531, ~l_5 = 33, ~l_6 = 4806021.$$
The corresponding higher elliptic unit $u_{1, \ok} = u_{2, \ok}^{-1} = \frac{v_2v_4v_6}{v_1v_3v_5}$ is close to the root $\approx -11.6360077... + i\cdot3.4634701...$ of the polynomial 
$$x^{10} + 24x^9 + 164x^8 + 99x^7 - 62x^6 - 89x^5 - 62x^4 + 99x^3 + 164x^2 + 24x + 1$$
which defines an absolute equation of $ \kk^{+}(\goth{f})$ over $\qq$. The computation time for $1000$ digits is 1 minute and 35 seconds, but the computation time for each of the individual computations is not uniform. The fifth computation requires 1 second whereas the sixth computation requires 58 seconds because of the level difference $l_5 = 33$ versus $l_6 = 4806021$.
\bigskip

\begin{table}[H]
\caption{Parameters for the quintic example}
\begin{center}
\noindent \begin{tabular}{| C C C C C C |}
\hline
\tau_1 \hfill= &935z^4 & - 19z^3 & - 2927z^2 & - 601z & + 796987 \\
\sigma_1 \hfill= &-3z^4 & + 1556z^3 & - 1205z^2 & - 3072z & + 987058 \\
\rho_1 \hfill = &3978z^4 & - 5242z^3 & + 1095z^2 & - 7073z & + 590241 \\
\varpi_1 \hfill = &-2767z^4 & + 4003z^3 & + 389z^2 & + 5232z & - 2072505\\
\hline

\tau_2 \hfill = &-377z^4 &+ 417z^3 & + 141z^2 & + 556z & - 19860 \\
\sigma_2 \hfill = &-71z^4 & - 74z^3 & + 105z^2 & + 449z & + 71171 \\
\rho_2 \hfill = &481z^4 & - 702z^3 & + 330z^2 & - 936z & - 324847 \\
\varpi_2 \hfill = &37z^4 & - 54z^3 & - 101z^2 & - 72z & + 28978\\
\hline

\tau_4 \hfill = &-121z^4 & + 180z^3 & + 17z^2 & + 254z & - 12881 \\
\sigma_4 \hfill = &830z^4 & - 1255z^3 & - 214z^2 & - 1580z & + 326894 \\
\rho_4 \hfill = &-6z^4 & - 56z^3 & + 82z^2 & + 41z & - 52774 \\
\varpi_4 \hfill = &-2082z^4 & + 3154z^3 & + 467z^2 & + 3916z & - 773076\\
\hline

\tau_3 \hfill = &-42z^4 & + 57z^3 & + 14z^2 & + 79z & - 139 \\
\sigma_3 \hfill = &287z^4 & - 443z^3 & - 60z^2 & - 522z & + 2216 \\
\rho_3 \hfill = &-3z^4 & + 27z^3 &+ z^2 & - 2z & - 522 \\
\varpi_3 \hfill = &-685z^4 & + 1029z^3 & + 157z^2 & + 1291z & - 2774\\
\hline

\tau_5 \hfill = &z^4 & - z^3 & + z^2 & - 4z & + 455 \\
\sigma_5 \hfill = & & &-2z^2 & + 3z & - 574 \\
\rho_5 \hfill = &-2z^4 & + 3z^3 & - 4z^2 &+ 8z & - 1282 \\
\varpi_5 \hfill = & & & z^2 & - z & + 247\\
\hline

\tau_6 \hfill = &-935z^4 & - 1775z^3 & + 1242z^2 & + 4691z & - 9783958 \\
\sigma_6 \hfill = &-2269z^4 & + 1923z^3 & + 6285z^2 & + 169z & + 19885070 \\
\rho_6 \hfill = &7142z^4 & - 1239z^3 & + 949z^2 & - 9197z & + 13206973 \\
\varpi_6 \hfill = &-11437z^4 & + 14892z^3 & + 1641z^2 & + 20621z & - 7776348\\
\hline
\end{tabular} 
\label{tablequintic}
\end{center}
\end{table}

\stepcounter{example}
\subsection{Example \exc: a discussion on a degree 6 example}

In \cite{thesis} we presented an example where we compute a higher elliptic unit above a degree $6$ field in an optimal setting. It takes a total of 6 pages to write down the $120$ evaluation parameters alone, so we will only briefly discuss this example here. The setting is the number field $\kk = \qq(z)$ where $\zexc$ is the complex root of the polynomial $x^6 - 2x^5 - 3x^4 + 10x^3 + 3x^2 - 8x - 3$ lying in the upper half-plane. The number field $ \kk = \qq(\zexc)$ has class number $1$. We choose $\goth{f}$ to be the unique prime ideal of norm $q = 2$. The corresponding narrow ray class group is $\classgroupf \simeq \zsz{2}$ and it is isomorphic to the narrow Hilbert class group of $\kk$. If we fix the fundamental units
\begin{align*}
\eps_1 &= 3\zexc^5 - 8\zexc^4 - 3\zexc^3 + 31\zexc^2 - 15\zexc - 11 \\
\eps_2 &= (2\zexc^5 - 5\zexc^4 - \zexc^3 + 18\zexc^2 - 13\zexc - 2)/5\\
\eps_3 &= (-4\zexc^5 + 10\zexc^4 + 7\zexc^3 - 41\zexc^2 + 6\zexc + 29)/5\\
\eps_4 &= (2\zexc^5 - 5\zexc^4 - 16\zexc^3 + 8\zexc^2 + 22\zexc + 8)/5
\end{align*}
for $\opck = \opc{\goth{f}}$ and the ordering of $\goth{S}_4$ given by Pari/GP, then we may compute the contents $\tlambda_1 = \dots = \tlambda_{24} = 1$ and the associated overflows

\begin{tabular}{CLCL}
\ttt_1 = & 3 \cdot821\cdot85146905507 &
\ttt_2 = &31\cdot71380217\cdot5479992107 \\
\ttt_3 = &2593\cdot28232090533 &
\ttt_4 = &4793\cdot20161\cdot11384677\\
\ttt_5 = &37\cdot18731\cdot12207408718823 & 
\ttt_6 = &757\cdot6354278197 \\
\ttt_7 = &5\cdot311\cdot4219\cdot261707099 & 
\ttt_8 = &5\cdot775267\cdot10124654046373 \\
\ttt_9 = &5^2\cdot37\cdot1061\cdot1321\cdot6449 &
\ttt_{10} = &5\cdot 64781598487 \\
\ttt_{11} = &5\cdot 9967\cdot 23337159379 &
\ttt_{12} = &5\cdot 1444658023 \\
\ttt_{13} = &426868799166283769 &
\ttt_{14} = &4950440701232129 \\
\ttt_{15} = &1361\cdot7219\cdot388087223 &
\ttt_{16} = &20173\cdot 200544349 \\
\ttt_{17} = &103\cdot 11593\cdot919759 &
\ttt_{18} = &523\cdot31751222611 \\
\ttt_{19} = &5\cdot 1759763531590697 &
\ttt_{20} = &15629\cdot3287441 \\
\ttt_{21} = &5\cdot 56807\cdot 354098824061 &
\ttt_{22} = &5\cdot73\cdot79\cdot349\cdot73009 \\
\ttt_{23} = &3544640951 &
\ttt_{24} = &293\cdot11719\cdot13931 
\end{tabular}  \smallskip

We may choose $\goth{b} = \ok$ and $\goth{a}$ the unique prime ideal of norm $N = 29$ in $ \ok$. We then compute the $24$ quotients: 
$$v_j = \frac{G_4\left(\frac{1}{2}, \frac{\tau_j}{\ell_j}, \frac{\sigma_j}{\ell_j}, \frac{\rho_j}{\ell_j}, \frac{\varpi_j}{\ell_j}, \frac{\xi_j}{\ell_j}\right)^{29}}{G_4\left(\frac{29}{2}, \frac{29\tau_j}{\ell_j}, \frac{29\sigma_j}{\ell_j}, \frac{29\rho_j}{\ell_j}, \frac{29\varpi_j}{\ell_j}\right)} $$
where the parameters $\tau_j, \sigma_j, \rho_j, \varpi_j, \xi_j$ are given in \cite{thesis} and we have defined for ease of presentation the levels $\ell_j = 2\cdot29 \cdot \ttt_j$. The corresponding higher elliptic unit 
$$u_{1, \ok} = \frac{v_2v_3v_6v_7v_{10}v_{11}v_{14}v_{15}v_{18}v_{19}v_{22}v_{23}}{v_1v_4v_5v_8v_9v_{12}v_{13}v_{16}v_{17}v_{20}v_{21}v_{24}}$$
coincides up to 1000 digits of precision with the root $\approx -555050859076374984.5110063... - i\cdot 334493188695056032.1307346...$ of the reciprocal polynomial 
\begin{align*}
P_{\mathrm{abs}} = & \,x^{12} + 1 + 1110101718152749974(x^{11}+x) \\ 
&+ 419967149444808248584979504483010229(x^{10} + x^2) \\
& + 2090591105457346230355038086262355202(x^9 +x^3) \\ 
&+ 4535292963058947524812988357459161338(x^8 + x^4)\\ 
&+ 6020269972074463320492578065966460430(x^7 + x^5)\\
& + 6311227681584443751632555661724851277x^6
\end{align*}
which defines an absolute equation of $\classfieldf$ over $\qq$.

\bibliographystyle{alpha}
\bibliography{bibliographie}{}

@misc{BCG,
	author = {Bergeron, Nicolas and Charollois, Pierre and Garc\'ia, Luis E.},
	title = {Elliptic units for complex cubic fields},
	eprint={2311.04110},
	archivePrefix={arXiv},
	journal ={},
	year = {2023},
	howpublished = {https://arxiv.org/abs/2311.04110}
}

@article{CharolloisDarmon,
 author = {Charollois, Pierre and Darmon, Henri},
 title = {Arguments des unit\'es de {S}tark et p\'eriodes de s\'eries d'{E}isenstein},
 fjournal = {Algebra \& Number Theory},
 journal = {Algebra Number Theory},
 issn = {1937-0652},
 volume = {2},
 number = {6},
 pages = {655--688},
 year = {2008},
 language = {French},
 doi = {10.2140/ant.2008.2.655},
 zbMATH = {5529341},
 Zbl = {1206.11062}
}

@article{DasguptaShintani,
 author = {Dasgupta, Samit},
 title = {Shintani zeta functions and {Gross}-{Stark} units for totally real fields},
 fjournal = {Duke Mathematical Journal},
 journal = {Duke Math. J.},
 issn = {0012-7094},
 volume = {143},
 number = {2},
 pages = {225--279},
 year = {2008},
 language = {English},
 doi = {10.1215/00127094-2008-019},
 url = {semanticscholar.org/paper/5e9c9b9b848dc95c431087454771b7fb11de9471},
 zbMATH = {5291473},
 Zbl = {1235.11102}
}

@article{CD,
 author = {Charollois, Pierre and Dasgupta, Samit},
 title = {Integral {Eisenstein} cocycles on {{\(\mathrm{GL}_n\)}}. {I}: {Sczech}'s cocycle and {{\(p\)}}-adic {{\(L\)}}-functions of totally real fields},
 fjournal = {Cambridge Journal of Mathematics},
 journal = {Camb. J. Math.},
 issn = {2168-0930},
 volume = {2},
 number = {1},
 pages = {49--90},
 year = {2014},
 language = {English},
 doi = {10.4310/CJM.2014.v2.n1.a2},
 zbMATH = {6324777},
 Zbl = {1353.11074}
}

@article{Ruijsenaars,
 author = {Ruijsenaars, S. N. M.},
 title = {First order analytic difference equations and integrable quantum systems},
 fjournal = {Journal of Mathematical Physics},
 journal = {J. Math. Phys.},
 issn = {0022-2488},
 volume = {38},
 number = {2},
 pages = {1069--1146},
 year = {1997},
 language = {English},
 doi = {10.1063/1.531809},
 url = {ir.cwi.nl/pub/2164},
 zbMATH = {1014649},
 Zbl = {0877.39002}
}

@article{FV,
 author = {Felder, Giovanni and Varchenko, Alexander},
 title = {The elliptic gamma function and {{\(\text{SL}(3,\mathbb Z)\ltimes\mathbb Z^3\)}}.},
 fjournal = {Advances in Mathematics},
 journal = {Adv. Math.},
 issn = {0001-8708},
 volume = {156},
 number = {1},
 pages = {44--76},
 year = {2000},
 language = {English},
 doi = {10.1006/aima.2000.1951},
 zbMATH = {1579912},
 Zbl = {1038.11029}
}

@article{FDuke,
 author = {Felder, Giovanni and Henriques, Andr{\'e} and Rossi, Carlo A. and Zhu, Chenchang},
 title = {A gerbe for the elliptic gamma function},
 fjournal = {Duke Mathematical Journal},
 journal = {Duke Math. J.},
 issn = {0012-7094},
 volume = {141},
 number = {1},
 pages = {1--74},
 year = {2008},
 language = {English},
 doi = {10.1215/S0012-7094-08-14111-0},
 url = {ora.ox.ac.uk/objects/uuid:8c728bbe-3e8c-4b4d-8c28-d68c5220a735},
 zbMATH = {5224873},
 Zbl = {1130.33010}
}

@article{Nishizawa,
 author = {Nishizawa, Michitomo},
 title = {An elliptic analogue of the multiple gamma function},
 fjournal = {Journal of Physics A: Mathematical and General},
 journal = {J. Phys. A, Math. Gen.},
 issn = {0305-4470},
 volume = {34},
 number = {36},
 pages = {7411--7421},
 year = {2001},
 language = {English},
 doi = {10.1088/0305-4470/34/36/320},
 zbMATH = {1695708},
 Zbl = {0993.33016}
}

@article{Narukawa,
 author = {Narukawa, Atsushi},
 title = {The modular properties and the integral representations of the multiple elliptic gamma functions},
 fjournal = {Advances in Mathematics},
 journal = {Adv. Math.},
 issn = {0001-8708},
 volume = {189},
 number = {2},
 pages = {247--267},
 year = {2004},
 language = {English},
 doi = {10.1016/j.aim.2003.11.009},
 zbMATH = {2136504},
 Zbl = {1077.33024}
}

@article{Sczech,
 author = {Sczech, Robert},
 title = {Eisenstein group cocycles for {{\(\text{GL}_ n\)}} and values of {{\(L\)}}- functions},
 fjournal = {Inventiones Mathematicae},
 journal = {Invent. Math.},
 issn = {0020-9910},
 volume = {113},
 number = {3},
 pages = {581--616},
 year = {1993},
 language = {English},
 doi = {10.1007/BF01244319},
 url = {https://eudml.org/doc/144138},
 zbMATH = {549667},
 Zbl = {0809.11029}
}

@manual{parigp,
	organization = "{The PARI~Group}",
	title = "{PARI/GP version \texttt{2.15.4}}",
	year = 2023,
	address = "Univ. Bordeaux",
	note = "available from http://pari.math.u-bordeaux.fr/"
}

@misc{lmfdb,
  shorthand    = {LMFDB},
  author       = {The {LMFDB Collaboration}},
  title        = {The {L}-functions and modular forms database},
  howpublished = {https://www.lmfdb.org},
  note         = {[Online; accessed 10 December 2025]},
}

@misc{firstpaper,
 author = {Morain, Pierre L. L.},
 title = {Geometric families of multiple elliptic {Gamma} functions and arithmetic applications, {I}},
 year = {2025},
 howpublished = {Preprint, {arXiv}:2510.16515 [math.{NT}]},
 url = {https://arxiv.org/abs/2510.16515},
 arXiv = {arXiv:2510.16515}
}

@misc{secondpaper,
 author = {Morain, Pierre L. L.},
 title = {Geometric families of multiple elliptic {Gamma} functions and arithmetic applications, {II}},
 year = {2026},
 howpublished = {Preprint, {arXiv}:2602.06561 [math.{NT}]},
 url = {https://arxiv.org/abs/2602.06561},
 arXiv = {arXiv:2602.06561}
}

@misc{thirdpaper,
	author = {Morain, Pierre L. L.},
	title = {Geometric families of multiple elliptic {Gamma} functions and arithmetic applications, {III} [{In} progress]},
	year = {2026}
}

@book{Neukirch,
	title = {Algebraic Number Theory},
	series = {Grundlehren der mathematischen Wissenschaften},
	author = {J. Neukirch},
	volume = {322},
	year = {1999},
	publisher = {Springer Berlin, Heidelberg},
	edition = {First}
}

@article{Robert,
 author = {Robert, Gilles},
 title = {Unit{\'e}s elliptiques et formules pour le nombre de classes des extensions ab{\'e}liennes d'un corps quadratique imaginaire},
 fjournal = {Bulletin de la Soci{\'e}t{\'e} Math{\'e}matique de France. Suppl{\'e}ment. M{\'e}moires},
 journal = {Bull. Soc. Math. Fr., Suppl., M{\'e}m.},
 issn = {0583-8665},
 volume = {36},
 pages = {77},
 year = {1973},
 language = {French},
 url = {https://eudml.org/doc/94657},
 zbMATH = {3491076},
 Zbl = {0314.12006}
}

@book{Cohen,
 author = {Cohen, Henri},
 title = {A course in computational algebraic number theory},
 fseries = {Graduate Texts in Mathematics},
 series = {Grad. Texts Math.},
 issn = {0072-5285},
 volume = {138},
 isbn = {3-540-55640-0},
 year = {1993},
 publisher = {Berlin: Springer-Verlag},
 language = {English},
 zbMATH = {435565},
 Zbl = {0786.11071}
}

@misc{rankoneStark,
	title = {Lecture notes on The {R}ank {O}ne {A}belian {S}tark {C}onjecture, {A}rizona {W}inter {S}chool},
	author = {S. Dasgupta and M. Greenberg},
	year = {2011},
}

@article{Espinoza,
	title = {Signed Shintani cones for number fields with one complex place},
	journal = {J. Number Theory},
	volume = {145},
	pages ={496-539},
	author = {M. Espinoza},
	year = {2014}
}

@article{roblot,
 author = {Roblot, Xavier-Fran{\c{c}}ois},
 title = {Stark's conjectures and {Hilbert}'s twelfth problem},
 fjournal = {Experimental Mathematics},
 journal = {Exp. Math.},
 issn = {1058-6458},
 volume = {9},
 number = {2},
 pages = {251--260},
 year = {2000},
 language = {English},
 doi = {10.1080/10586458.2000.10504650},
 url = {https://eudml.org/doc/224431},
 zbMATH = {1647647},
 Zbl = {0986.11074}
}

@article{Ramachandra,
 author = {Ramachandra, K},
 title = {Some applications of {Kronecker}'s limit formulas},
 fjournal = {Annals of Mathematics. Second Series},
 journal = {Ann. Math. (2)},
 issn = {0003-486X},
 volume = {80},
 pages = {104--148},
 year = {1964},
 language = {English},
 doi = {10.2307/1970494},
 zbMATH = {3229683},
 Zbl = {0142.29804}
}

@phdthesis{thesis,
  TITLE = {{Higher elliptic Gamma functions, cohomology of linear groups and higher elliptic units}},
  AUTHOR = {Morain, Pierre L. L.},
  URL = {https://hal.science/tel-05674839},
  NUMBER = {2026SORUS092},
  SCHOOL = {{Sorbonne Universit{\'e}}},
  YEAR = {2026},
  TYPE = {Thesis},
  HAL_ID = {tel-05674839},
  HAL_VERSION = {v2},
  note = {Available on HAL at https://hal.science/tel-05674839} 
}

@article{CassouNogues,
 author = {Cassou-Nogu{\`e}s, Pierrette},
 title = {Values at negative integers of zeta functions and {{\(p\)}}-adic zeta functions},
 fjournal = {Inventiones Mathematicae},
 journal = {Invent. Math.},
 issn = {0020-9910},
 volume = {51},
 pages = {29--59},
 year = {1979},
 language = {French},
 doi = {10.1007/BF01389911},
 url = {https://eudml.org/doc/142621},
 zbMATH = {3634372},
 Zbl = {0408.12015}
}

\end{document}